\documentclass[12pt]{amsart}

\usepackage{amssymb}
\usepackage{amsmath}
\usepackage{amsthm}
\usepackage{amssymb}
\usepackage{marginnote}
\theoremstyle{plain}

\newtheorem*{theorem*}{Theorem}

\newtheorem*{proposition*}{Proposition}
\newtheorem{lemma}[subsubsection]{Lemma}
\newtheorem*{lemma*}{Lemma}

\newtheorem*{corollary*}{Corollary}

\newtheorem*{claim}{Claim}
\theoremstyle{definition}
\newtheorem*{definition}{Definition}
\theoremstyle{remark}
\newtheorem*{remark}{Remark}
\newtheorem*{remarks}{Remarks}
\newtheorem{remarknum}[subsubsection]{Remark}

\newtheorem*{eexample*}{Example}

\newcommand{\OO}{\mathcal{O}}

\newcommand{\Hom}{\operatorname{Hom}}
\newcommand{\End}{\operatorname{End}}

\newcommand{\GL}{\operatorname{GL}}
\newcommand{\GLt}{\operatorname{GL}(2)}

\newcommand{\SL}{\operatorname{SL}}

\newcommand{\tr}{\operatorname{tr}}

\newcommand{\ev}{\operatorname{ev}}
\newcommand{\meas}{\operatorname{\meas}}
\newcommand{\stm}{\setminus}
\newcommand{\nc}{\newcommand}
\nc{\cal}{\mathcal} 
\nc{\la}{\langle} \nc{\ra}{\rangle}
 \nc{\CA}{\cal A}
 \nc{\CBB}{\cal B}
\nc{\CDD}{\cal D}
\nc{\CE}{\cal E}
\nc{\CF}{\cal F} \nc{\CG}{\cal
G} \nc{\CH}{\cal H} \nc{\CI}{\cal I} \nc{\CJ}{\cal J}
\nc{\CK}{\cal K} \nc{\CL}{\cal L} \nc{\CM}{\cal M} \nc{\CN}{\cal
N} \nc{\CO}{\cal O} \nc{\CP}{\cal P} \nc{\CQ}{\cal Q}
\nc{\CR}{\cal R} \nc{\CS}{\cal S} \nc{\CT}{\cal T} \nc{\CU}{\cal
U} \nc{\CV}{\cal V} \nc{\CW}{\cal W} \nc{\CZ}{\cal Z}

\nc{\Ck}{\textsl{k}}

\nc{\fa}{\mathfrak a} \nc{\fg}{\mathfrak g} \nc{\fii}{\mathfrak i}\nc{\fk}{\mathfrak k}
\nc{\fh}{\mathfrak h} \nc{\fm}{\mathfrak m} \nc{\fn}{\mathfrak n}  \nc{\fG}{\mathfrak G}
\nc{\fA}{\mathfrak A} \nc{\fC}{\mathfrak C} \nc{\fI}{\mathfrak I}
\nc{\fL}{\mathfrak L} \nc{\fS}{\mathfrak S}
\nc{\fz}{\mathfrak z} \nc{\fl}{\mathfrak l}
\nc{\fp}{\mathfrak p}
\nc{\ft}{\mathfrak t}

\nc{\nen}{\newenvironment} \nc{\ol}{\overline}
\nc{\ul}{\underline} \nc{\lra}{\longrightarrow}
\nc{\lla}{\longleftarrow} \nc{\Lra}{\Longrightarrow}
\nc{\Lla}{\Longleftarrow} \nc{\Llra}{\Longleftrightarrow}
\nc{\hra}{\hookrightarrow} \nc{\iso}{\overset{\sim}{\lra}}

\makeatletter
\@addtoreset{equation}{section}
\makeatother
\numberwithin{equation}{section}

 \nc{\ba}{\mathbb A}
 \nc{\bq}{\mathbb Q}
 \nc{\br}{\mathbb R}
 \nc{\bz}{\mathbb Z}
 \nc{\bc}{\mathbb C}
 \nc{\bn}{\mathbb N}
\nc{\bg}{\mathbb G}
 \nc{\ck}{\mathcal{K}}
 \nc{\G}{\Gamma}
 \nc{\sm}{\setminus}
 \nc{\sub}{\subset}
 \nc{\lm}{\lambda}
  \nc{\Lm}{\Lambda}
 \nc{\al}{\alpha}
 \nc{\bt}{\beta}
 \nc{\om}{\omega}
 \nc{\dl}{\delta}
 \nc{\g}{\gamma}
 \nc{\Dl}{\Delta}
 \nc{\Om}{\Omega}
 \nc{\s}{\sigma}
 \nc{\ro}{\rho}
 \nc{\te}{\theta}
 \nc{\SLR}{\operatorname{SL}_2(\br)}
 \nc{\GLR}{\operatorname{GL}_2(\br)}
 \nc{\PGLR}{\operatorname{PGL}_2(\br)}
 \nc{\PSLR}{\operatorname{PSL}_2(\br)}
 \nc{\PSLZ}{\operatorname{PSL}_2(\bz)}
 \nc{\SLC}{\operatorname{SL}_2(\bc)}
 \nc{\uH}{\mathbb H}
 \nc{\fD}{\mathcal{D}}
 \nc{\fE}{\mathcal{E}}
 \nc{\fO}{\mathcal{O}}
 \nc{\haf}{\frac{1}{2}}
 \nc{\qtr}{\frac{1}{4}}
 \nc{\shaf}{{\scriptstyle\frac{1}{2}}}
 \nc{\hlm}{{\scriptstyle\frac{\lambda}{2}}}

 \nc{\8}{\infty}
 \nc{\7}{{-\infty}}
 \nc{\inv}{^{-1}}
 \nc{\eps}{\varepsilon}
 \nc{\aG}{\mathbf{G}}
 \nc{\spn}{\operatorname{Span}}
 \nc{\Cm}{\operatorname{CM}}
 \nc{\tildl}{\dl^1}
\nc{\chiv}{{\chi_\fp}}
\nc{\psiv}{{\psi_\fp}}
\nc{\piv}{{\pi_\fp}}
\nc{\zt}{Z\setminus T}

\newcommand*\MapsTo{%
  {\xrightarrow[\raisebox{0.25 em}{\smash{\ensuremath{\sim}}}]{}}%
}

\begin{document}
\title[Periods and Invariants]{Periods and global invariants of automorphic representations}

\author{Joseph Bernstein} \email{}

\address{School of Mathematics\\
Tel Aviv University \\
Ramat Aviv, 69978\\
Israel}

\author{Andre Reznikov} \email{}

\address{Department of Mathematics\\
Bar-Ilan University \\
Ramat Gan, 52900 \\
Israel}

\begin{abstract}We consider periods of automorphic representations of adele groups  defined by integrals along Gelfand subgroups. We define natural maps between local components of such periods and construct corresponding  global maps using automorphic $L$-functions. This leads to an introduction of a global invariant of an automorphic representation arising from two such periods. We compute this invariant in some cases.

\end{abstract}


\thanks{The research  was partially supported  by the ERC grant  291612, by the ISF   grant 533/14, and, during the visit to IAS,  by the National Science Foundation under Grant No. DMS - 1638352.}
\maketitle

\section{Introduction}\label{intro-sect} 
\subsection{Periods and  relations between them.} Periods play a central role in the modern theory of automorphic functions. In particular, there are many instances where periods of automorphic functions are related to $L$-functions. In this paper, we propose to study relations between different periods defined on the same  automorphic representation. We find that in this new setup $L$-functions also appear naturally. 
 
 Our main aim in this paper is to highlight a relation between different periods which, we believe underlies ``period to $L$-function'' relation (in fact, the celebrated paper of J-L. Walsdpurger \cite{Wa} that was our starting point,  already contains the idea we are trying to expand). By doing so, we are able make a very natural reformulation of of some classical results and to include in our scheme periods which, apparently, do not fall in the familiar framework  (see Section \ref{opp-whittak-sect}).

We explore a general scheme that allows to compare different periods and consider few (classical and not so classical) examples. We first consider relations between  the Whittaker (i.e., unipotent) period and the Hecke (i.e., torus) period for $\GLt$. One of the relations (see Theorem \ref{N-per-T-per-thm}) is classical and is a mere reformulation of the treatment given by  H. Jacquet and R. Langlands \cite{JL} to the Hecke method. However, we discover a converse relation (see Theorem \ref{T-per-N-per-thm}) which seems to be new (although similar local considerations appeared recently in \cite{SV}). We then consider a non-classical example of two Whittaker periods on different unipotent subgroups of $\GLt$.  In this case, our construction leads to an Euler product with a non-standard local factor which nevertheless could be regularized with the help of an appropriate $L$-function (see Section \ref{opp-whittak-sect}).  This leads to an introduction of a non-trivial {\it global invariant} of an automorphic representation of $\GLt$. 

Our motivation comes from the desire to understand the fundamental ``period to $L$-function'' relation in the theory of automorphic functions. 
Such a relation goes back to the foundational work of E. Hecke \cite{He}, where he constructed the Hecke $L$-function on $\GLt$ as the period integral along the split torus in $\GLt$.  This is the most basic of  ``period to $L$-function'' relations. Another striking example was discovered by J.-L. Waldspurger \cite{Wa} and  connects the period along a non-split torus in $\GLt$ to the special value of an $L$-function of the appropriate base change lift (i.e., a period on another group). Our starting point was an attempt to reformulate Waldspurger's theorem in representation theoretic terms which do not require the notion of $L$-function.  

 We mention a vast generalization of the Waldspurger's result formulated as a conjecture by B. Gross and D. Prasad \cite{GP}.      Consequently, the exact form of the Gross-Prasad period relation  was conjectured by A. Ichino and T. Ikeda \cite{II}. This led to other formulas relating normalized periods and $L$-functions (e.g., an analog for the Whittaker functional was considered in \cite{LM1}).  A general framework for period formulas in the context of Plancherel measures was recently proposed by Y. Sakellaridis and A. Venkatesh \cite{SV}. We also note recent work \cite{LM2} where similar to ours ideas discussed for much more sophisticated higher rank examples (under the name of ``model transition''). However, our interest is slightly different from \cite{LM2} as  we are interested in understanding the nature of the global invariant associated to two periods.

We note that one of the most important attributes of period to $L$-function formulas is the presence of the multiplicity one phenomenon (i.e., the Gelfand property of one-dimensionality of certain invariant functionals; see   \cite{Gr}). This point of view was pioneered by I. Piatetski-Shapiro \cite{PS}, and also will be  essential throughout this paper.  
\subsection{Action on periods} We are interested in the following setup. Let $G$ be an algebraic (reductive) group over a global field $k$ (in practice a reader can assume $k=\bq$ for simplicity), and  let $H_1, H_2 \subset G$  be  two algebraic subgroups of $G$ also defined over $k$ (e.g., a split over $k$ torus and an associated unipotent subgroup in $G=\GLt$).  Let $G(\ba), H_1(\ba), H_2(\ba)$ be the corresponding adele groups, and we denote by $X_G=G(k)\setminus G(\ba), X_{H_1}={H_1}(k)\setminus {H_1}(\ba), X_{H_2}={H_2}(k)\setminus {H_2}(\ba)$ the corresponding automorphic quotient spaces.  Let $\pi$ be an automorphic representation of $G$ (we will be vague at this point of what is required of $\pi$). 
We are interested in the period functional given by the integral $p_{H_1}(\phi)=\int _{X_{H_1}}\phi(h)dh$ over the ${H_1}(\ba)$-orbit  $X_{H_1}\subset X_G$ of an automorphic function $\phi$ belonging to the space of the representation  $\pi$ (and similarly for the period $p_{H_2}$ for $X_{H_2}$). More generally, we consider periods twisted by  characters $\chi_i:{H_i}(k)\setminus {H_i}(\ba)\to\bc$ which are given by integrals $p_{H_i,\chi_i}(\phi)=\int _{X_{H_i}}\chi_i\inv(h_i)\phi(h_i)dh_i$. To define such periods, one has to choose (invariant) measures on subgroups and impose certain restrictions on representation $\pi$ and on spaces $X_{{H_i}}$. 

Assuming that all these periods are well-defined, it is natural to ask if there is a relation between functionals $p_{H_1}$ and $p_{H_2}$ which are defined on the same space of smooth vectors for the automorphic representation $\pi$. Periods $p_{H_1}$ and $p_{H_2}$ define functionals on $\pi$, and one possibility would be to compute their correlation (i.e., the scalar product, properly understood of course). In fact it is possible in many cases (see \cite{Gr}), but we found it a little bit easier to  make another comparison in terms of the action of adelic groups. Namely, we can try to integrate 
the functional $p_{H_1}$ with respect to the action of the adelic group ${H_2}(\ba)$.
 Assuming that such an operation is well-defined, we would obtain an ${H_2}(\ba)$-invariant functional $\tilde p_{H_2}=\int_{h\in _{H_2}(\ba)}\pi^*(h)p_{H_1}\ dh$ on $\pi$ (i.e., $\tilde p_{H_2}(v)=\int_{h\in _{H_2}(\ba)}\int_{x\in X_{H_1}}v(xh)dxdh$ for any smooth vector $v$ in the representation $\pi$). This does not identify such a functional in general, but in the case when ${H_2}(\ba)$ is a Gelfand subgroup of $G(\ba)$ (i.e., the space of ${H_2}(\ba)$-invariant functionals on $\pi$ is at most one-dimensional), we should get a functional which is proportional to the period functional $p_{H_2}$.  What we found is that the above mentioned ``classical'' period to $L$-function formulas allows one to compute the coefficient of proportionality between $\tilde p_{H_2}$ and $p_{H_2}$ in some cases. Moreover, we find  the ``$L$-functions free'' formulation of this relation between periods even more interesting. Such a reformulation allows us to consider cases where the relation to $L$-functions is not a standard one.


\subsection{The construction}\label{idea} We will work only with periods satisfying the local uniqueness property (and hence also satisfying global uniqueness). Let ${H}\subset G$ be a subgroup of a group both defined over a global field $k$. For a place $\fp$ of $k$, we consider local groups ${H_{\fp}}\subset G_\fp$,  (i.e.,  groups of points over the local field $k_\fp$). Let   $\pi=\hat\otimes\pi_\fp$ be an irreducible representation of $G(\ba)$ and $\chi=\hat\otimes\chi_\fp$ be a character of ${H}(\ba)$ (more generally, one can consider an irreducible representation of $H(\ba)$ as well). We consider the complex vector space of equivariant maps, {\it the periods space},  $\CP(\pi,\chi)=\Hom_{{H}(\ba)}(\pi,\bc_\chi)$ and its local counterparts, {\it the local period space}, $\CP_\fp(\pi_\fp,\chi_\fp)=\Hom_{{H_{\fp}}}(\pi_\fp,\bc_{\chi_\fp})$.
We call a tuple $(G_\fp,\pi_\fp,H_\fp,\chi_\fp)$ a local Gelfand  data (or a multiplicity one tuple) if $\dim \CP_\fp(\pi_\fp,\chi_\fp)\leq 1$. In such a case we have $\CP(\pi,\chi)=\hat\otimes \CP_\fp(\pi_\fp,\chi_\fp)$, and the global period space  is also at most one-dimensional. We call the tuple $(G,\pi,H,\chi)$ globally Gelfand if it is locally Gelfand at every place. We note that in practice we consider a slightly different space of $H$-maps from $\pi$ with values in co-invariants of $H$ (see Section \ref{periods-setup}). We find the language of co-invariants more appropriate when dealing with periods, and it leads to more canonical constructions.

Let $(G,\pi,{H_1},\chi_1)$
and $(G,\pi,{H_2},\chi_2)$ be two globally Gelfand tuples and $P(\pi,\chi_i)$ corresponding {\it one-dimensional} complex vector spaces. 
Our goal is to construct a {\it canonical} map 
\begin{eqnarray}\label{I}
I: \CP(\pi,\chi_1)\to \CP(\pi,\chi_2)
\end{eqnarray}
between these one-dimensional vector spaces in the presence of the corresponding  {\it automorphic} periods. We do this in two steps. 

\subsection{Local step}\label{local-step} The first step is purely local. It is relatively easy to construct {\it local} maps $I_\fp: \CP_\fp(\pi_\fp,\chi_{1,\fp})\to \CP_\fp(\pi_\fp,\chi_{2,\fp})$ between local spaces of periods using the integration along the subgroup $H_{2,\fp}\subset G_\fp$. For a given vector $\xi_\fp\in \CP_\fp(\pi_\fp,\chi_{1,\fp})$, we define a vector
$I_\fp(\xi_\fp)\in \CP_\fp(\pi_\fp,\chi_{2,\fp} )$ by $I_\fp(\xi_\fp):=\int_{{H_2,}_\fp}\chi_{2,\fp}\inv(h_\fp)\pi_\fp^*(h_\fp)(\xi_\fp)dh_\fp$, where $\pi^*_\fp$ denotes the dual representation of $G_\fp$ on $V^*_{\piv}$. The integral is understood in a weak sense. This means that for any smooth vector $v_\fp\in V_{\piv}$, we have
$I_\fp(\xi_\fp)(v_\fp)=\int_{H_{2,\fp}}\chi_{2,\fp}\inv(h_\fp)\xi_\fp(\pi_\fp(h_\fp)v_\fp)dh_\fp$ (in a case the intersection $H_{1,\fp}\cap H_{2,\fp}$ in non-trivial one have to take integration over the appropriate quotient). The last integral might be divergent, but in many cases could be  evaluated by a standard procedure (see \cite{G1} for a general strategy based on analytic continuation and \cite{LM3} for a compact exhaustion ubiquitous in a non-archimedian setup). 
We stress  that local maps are assumed to be defined {\it canonically} for {\it all} $\ \fp$. This of course might constitute a non-trivial local problem, but it is solved in many cases.

\begin{remark} We note that the alluded local map $I_\fp$ could be constructed (or normalized) in various ways. In particular, a very general method was envisioned by H. Jaqcuet by introducing the Relative Trace Formula (see \cite{Jacq-RTF} for the general framework). We hope to discuss this connection elsewhere.
\end{remark}

\subsection{Global step} The next step is to ``glue'' local maps $I_\fp$ to form a global map $I$. This is a more subtle procedure. 
We construct the global map $I$ by regularizing the tensor product $\otimes I_\fp$ of local maps with the help of appropriate weight factors. This is possible only for local maps which are coming from  automorphic periods, and the weight factors are provided by the theory of automorphic $L$-functions. 
The construction of the map $I$ (in certain cases) is the main observation of the paper.

\subsubsection{Regularization by $L$-functions}\label{reg-by-L} We propose the following procedure, having its origin in \cite{Wa}, for the regularization of the map  $\otimes I_\fp$. 
Let $\xi$ be a vector in $\CP(\pi,\chi_1)$. We  write it as a product $\xi=\otimes_\fp \xi_\fp$, where
$\xi_\fp\in \CP_\fp(\piv, \chi_{1,\fp})$ and for almost all $\fp$ we have $\xi_\fp(e^0_\fp) = 1$ for the standard vector $e_\fp^0\in V_\piv$ (see Section \ref{loc-map-coinv} for a slightly more precise formulation in terms of co-invariants). Given a decomposable vector $v\in V_\pi$, we  write it in a form $v = \otimes_\fp v_\fp$, where $v_\fp = e^0_\fp$ for almost all $\fp$.
Now we would like to set $I(v) =\prod_\fp d_\fp(v_\fp)$, where $d_\fp(v_\fp) := I_\fp(\xi_\fp)(v_\fp)$.

This product is usually not convergent. 
To regularize such an infinite product we can use the fact that outside
of a finite number of places maps $d_\fp$ can be explicitly computed using an unramified computation.  Unramified factors  do not  depend on a choice of the vector $v$. In many cases maps $d_\fp$ are connected to local $L$-factors of automorphic representations. In such cases this allows us to use the following regularization procedures.

The regularization procedure is based on a renormalization of  local unramified maps $d_\fp$ by scalar factors which we call $L$-weights. These $L$-weights are constructed from Langlands $L$-factors. 

First, assume we can find an "$L$-weight" consisting of an $L$-function, i.e., an appropriate (partial) automorphic $L$-function given by  $L_S(s) =
\prod_{\fp\not\in S}
L_\fp(s)$ for   $Re(s)\gg 1$,  and a complex number $s_0$ such that $d_\fp(e^0_\fp)= L_\fp(s_0)$ for $\fp\not\in S$. Here $S$ is a finite set of primes and for a fixed $\piv$, the Euler factors
for all $\fp\not\in S$ are  rational functions of $q_\fp^{-s}$ given by a local Langlands $L$-function.   Note that an existence of the factor $L_\fp(s_0)$ as above is a condition on the unramified representation $\pi_\fp$ and does not depend on the vector $v_\fp$.   We define the normalized maps $d^0_\fp(v_\fp): = L_\fp(s_0)\inv d_\fp(v_\fp) $. In particular, for an unramified $\fp$, we have  $d^0_\fp(e^0_\fp)=1$. Here we assume that $L_\fp(s_0)\not =0$ for all but a finite number of primes $\fp$.
We define now for a large enough finite set of  primes $S$,
\begin{eqnarray}\label{L-reg-intro-0}I(v) := L_S(s_0)\prod_{\fp\in S} d_\fp(v_\fp)\ .
\end{eqnarray} 
Here the set $S$ depends on the vector $v=\otimes_\fp v_\fp$.
In many cases, the complex number $s_0$ belongs to the region of the analytic continuation of $L(s)$, an extension to which we will take for granted.

\begin{remark} More generally, one have to  consider $L$-weights given by a ratio of various Langlands $L$-functions evaluated at possibly different points instead of just one $L$-function as above. We defer to Section \ref{dual-gr-fn} a discussion of why one might expect to see such weights in a regularization procedure. 
\end{remark}
The regularization \eqref{L-reg-intro-0} is in fact the standard way to define infinite Euler products in many instances arising in the theory of automorphic functions.  Our new point is that it could be extended to cover more cases which are less standard for the theory of automorphic functions (although the basic argument is well-known in other situations in number theory). 

Namely, we propose a slight elaboration of the above method. In Section \ref{opp-whittak-sect}  we will provide an example where  a local $L$-function (or a ratio of $L$-functions) such that  $d_\fp(e^0_\fp)= L_\fp(s_0)$ for $\fp\not\in S$ does not exists. On the other hand,  
we can find a partial $L$-function (or an $L$-weight) and some complex number $s_0$ such that
if we replace for almost all $\fp$, maps $d_\fp$ by  $d^0_\fp: = L_\fp(s_0)\inv d_\fp $, then the
product $\prod_\fp d^0_\fp(v_\fp)$ is absolutely convergent (this condition does not depend on a specific
choice of the vector $v=\otimes_\fp v_\fp$). 

Hence we propose the following modified regularization procedure.   For a large enough finite set of  primes $S$ (depending on a vector $v$ ), we define
\begin{eqnarray}\label{L-reg-intro}I(v) := L_S(s_0)\prod_{\fp\not\in S} d_\fp^0(v_\fp)\prod_{\fp\in S} d_\fp(v_\fp)\ . 
\end{eqnarray}  Note that \eqref{L-reg-intro} is a generalization of \eqref{L-reg-intro-0} since in the former case we have $ d^0_\fp(v_\fp)=1$ for ${\fp\not \in S}$, and hence can omit the middle term in the product.    

 In principle, there might be different regularizing factors $L_S$ giving for the same periods a different regularized map. We expect however that this is not the case and in fact each period map could be regularized by essentially unique combination of $L$-functions (if a regularization exists at all which we do not claim in general).

The basic idea behind regularization \eqref{L-reg-intro} is well known since the work of  T. Estermann \cite{E} (see also \cite{K1}, \cite{K2}) and usually does not produce interesting Euler products. Our point here is that notwithstanding this there is a {\it natural} construction in the theory of automorphic functions  which leads to such Euler products.

\subsection{Period invariant} Having constructed the map $I$, we can ask what is its effect  on automorphic periods. Namely, we can try to compare the original period functional  $p_{H_2,\chi_2}$ and the newly constructed functional $\tilde{p}_{H_2,\chi_2}=I(p_{H_1,\chi_1})$. This is the last (and the most interesting) step of the construction.  The coefficient of proportionality (when defined) gives rise to a {\it global invariant} of the automorphic representation $\pi$ (for $\chi_1$ and $\chi_2$ fixed). When $\tilde p_{H_2,\chi_2}=p_{H_2,\chi_2}$, this invariant is equal to $1$, and we say that the collection $\{I_\fp\}$ of local maps (or the resulting global map $I$) is {\it coherent}, that is, the {\it local to global principle} is satisfied for these automorphic periods. One can show that in many classical examples this is indeed the case, and this is equivalent to the ``period to $L$-function'' relation we mentioned above (e.g., theorems of Hecke and Waldspurger on torus periods).
However, we find that sometimes the relation between $\tilde p_{H_2,\chi_2}$ and $p_{H_2,\chi_2}$ is more complicated and this gives rise to a non-trivial invariant. In particular, in  Section \ref{opp-whittak-sect} we provide an example of two opposite unipotent subgroups of $GL(2)$ and Whittaker functionals $\tilde p_{H_2,\chi_2}$ and $p_{H_2,\chi_2}$  which do {\it not} coincide for the  cusp form associated with the Ramanujan tau function.  As a result we obtain some new, as far as we know, numerical invariant canonically associated to the Ramanujan cusp form (see Appendix \ref{Ramanj-apx}).


\subsection{Structure of the paper} In this paper, we start with a discussion of examples related to the classical periods considered  by Hecke (and in the adelic setting by H. Jacquet and R. Langlands \cite{JL}). We will show (see Theorem \ref{N-per-T-per-thm}) how to define a procedure of integration transforming  the Whittaker period (i.e., the period defining a non-trivial Fourier coefficient along the horocycle) into the Hecke period (i.e., the period against a Hecke character along the split torus of $\GLt$). The proof we present is a simple reformulation of the standard argument of Hecke-Jacquet-Langlands and is based on the unfolding technique. 

Next we  prove the converse statement (i.e., we integrate the Hecke period  into the Whittaker period; see Theorem \ref{T-per-N-per-thm}). This case is not covered by the Hecke-Jacquet-Langlands theory and the result is new. Here our proof is purely local (combined with the direct statement for the Whittaker to Hecke transform), and does not involve an unfolding procedure.  In fact, we do not know if some kind of an unfolding exists in this case. 

In Section \ref{opp-whittak-sect} we consider two unipotent periods for $\GLt$, that is,  two Whittaker functionals: one for $N^+=\left\{
\left(\begin{smallmatrix}
1 & x\\
                         &    1      \\
\end{smallmatrix} \right)\right\}$ and another for $N^-=\left\{
\left(\begin{smallmatrix}
1 & \\
                     x    &    1      \\
\end{smallmatrix} \right)\right\}$. We consider the same question as before and discover  that this example is of a completely different nature than those we discussed so far. We define the local integration procedure and show how to regularize the global map with the help of the adjoint $L$-function $L(1,\piv,Ad)$.  In this case, the relation of the period map to special values of $L$-functions is more puzzling as compared to the Hecke case. We note that local coefficients $d_\fp$ appearing in this case do {\it not} coincide with some familiar Euler factors from the theory of $L$-functions; however, coefficients $d_\fp$ coincide with the {\it linear part}  of the Euler polynomial of $L(1,\piv,Ad)$ (or, as one might say, with the leading term  of $L(1,\piv,Ad)$). This will be essential for the regularization of the Euler product $\prod_\fp d_\fp$.  
 We will show that there is natural map between Whittaker periods on different unipotent subgroups, but in general it {\it does not} map the automorphic period to the automorphic period. As a result we are able to define a non-trivial invariant (a "defect" or an "obstruction") of an automorphic (cuspidal) representation of $\GLt$. We also present a numerical computation for this invariant for the Ramanujan cusp form which indicates that the resulting invariant is not trivial (see Appendix~\ref{Ram}).

We would like to point out that the above phenomenon seems to be  part of a pattern and not an isolated example. In several other
instances, we have computed analogous local maps and found that these are connected to $L$-functions in a similar way (that is, coincide with linear parts of some $L$-functions). This should allow one to define the corresponding global maps (e.g., a map between a  torus period and a {\it non-associated} Whittaker period). We will discuss these examples elsewhere. 

We discuss in Section \ref{aver-val-L} a geometric reason why the particular $L$-factor $L(1,\piv,Ad)$ might be relevant for the regularization procedure for two Whittaker functionals. 
In Section \ref{dual-gr-fn}  we discuss a reason why one might expect Langlands $L$-functions be relevant to the regularization procedure we propose.

 In Section \ref{integration},
we review the basic setup and in particular discuss machinery of co-invariants which we find to be convenient in our treatment of periods of automorphic representations. 

In Appendix \ref{Kir-mod-apx}, we collect information about the Kirillov model for $\GLt$ which we use in proofs and for computations. 

\subsubsection{Notations} We denote the global field by $k$, places of $k$ by $\fp$, the set of places of $k$ by $\CP(k)$, the corresponding ring of adeles by $\ba$, and the  group of ideles by $J_\ba$. For a group $\CG$ defined  over a global field $k$ (e.g., over $k=\bq$), we denote by $\CG(k)$ the group of $k$-points, by $\CG_\fp=\CG(k_\fp)$ the group of points over a local field $k_\fp$ (e.g., over $k_\fp=\bq_p$ or $k_\8=\br$) and by $\CG_\ba=\CG(\ba)$ the group of adelic points. For a local non-archimedian field $k_\fp$, we denote by $\CO_\fp$ the ring of integers, by $\varpi_\fp$ a generator of the maximal ideal in $\CO_\fp$, and by $q_\fp$ the size of the residue field. We will use the letter $\psi$ to denote  additive characters (local or global), and the letter $\chi$ to denote multiplicative characters (local or global). For a quasi-character $\chi_\fp:k_\fp^\times\to\bc^\times$ of a local field, we have the decomposition  $\chi_\fp=|\chi_\fp|\cdot \tilde\chi_\fp$ where $\tilde\chi_\fp$ is unitary.    We denote by $\s_{\chi_\fp}=Re(\chi_\fp)\in\br$ the real part of $\chi_\fp$ given by the relation $|\chi_\fp(\varpi_\fp)|=|\varpi_\fp|^{\s_{\chi_\fp}}$. 
Similarly, for a Hecke character $\chi:k^\times\stm J_\ba\to \bc^\times$, there exists the unique decomposition $\chi=|\chi|\cdot\tilde\chi$ and $|\chi|=|\cdot|^{\s_\chi}$ with  $\s_\chi=Re(\chi)\in\br$. 

We denote  by $G=\GLt$, by $Z=Z_G$ its center, by $T$ the subgroup of diagonal matrices, by $A\subset T$ the subgroup of matrices of the form $
\left(\begin{smallmatrix}
a & \\
                         &    1      \\
\end{smallmatrix} \right)$, and by $N=N^+$ the upper triangular matrices. We will  use  the following notations: $n(x)=
\left(\begin{smallmatrix}
1 & x\\
                         &    1      \\
\end{smallmatrix} \right)$,
$\bar a=
\left(\begin{smallmatrix}
a & \\
                         &    1      \\
\end{smallmatrix} \right)$,
$diag(a,b)=
\left(\begin{smallmatrix}
a & \\
                         &    b      \\
\end{smallmatrix} \right)$, 
$z(a)=
\left(\begin{smallmatrix}
a & \\
                         &    a      \\
\end{smallmatrix} \right)$,
$w=
\left(\begin{smallmatrix}
 & -1\\
                  1      &        \\
\end{smallmatrix} \right)$.   

For a field $F$ (e.g., $F=k_\fp$), characters $\chi:A(F)\to \bc^\times$ are given by $\chi(\bar a)=\chi(a)$, $a\in F$, and hence we can identify  characters of $A$ with those of $F^\times$. We do this for global characters as well and hence identify Hecke characters of $k$ with those of $A(\ba)$.  We use the notion of the real part for  local and global characters of $A$ as well.

\subsubsection{Acknowledgments}  We would like to thank E. Baruch, B. Gross, E.~Lapid, S. D. Miller, M. Rubenstein,  Y. Sakellaridis, D. Zagier and the referee for  valuable comments.


\section{Measures, automorphic representations and periods}\label{integration}

\subsection{Invariant integration} We review invariant measures on local and adelic groups. We make an emphasis on the language of co-invariants. 

\subsubsection{Torsors} By a torsor we mean a one-dimensional complex vector space. The name comes from the fact that if $L$ is a torsor then $L\setminus \{0\}$ is a $\bc^\times$-torsor. Torsors form a tensor category with respect to the tensor product. This category has the unit object $L_0=\bc$, and for every torsor $L$ there is an inverse torsor $L\inv:= L^*$. 

\subsubsection{Moderate groups} Let $\CA$ be a locally compact group. We say that $\CA$ is moderate if there exists a compact subgroup $K \subset \CA$ with the following properties:

\begin{enumerate}
	\item[(i)] $K$ is totally disconnected,
	\item[(ii)] The normalizer $N$ of the group $K$ is open in $\CA$, and the quotient group $N/K$ is a (smooth) Lie group.
\end{enumerate}
We call a subgroup $K$ with these properties a basic compact subgroup.

  We will work only with moderate groups. In fact, as follows from Gleason-Yamabe's theorem  (see \cite{T}, Exercise 1.6.4,  \cite{MZ}, p. 182), any
locally compact group of a finite topological dimension is moderate.

\begin{proposition*}Any two basic compact subgroups $K, L$ in a moderate group $\CA$ are commensurable, i.e.,
the group $L\bigcap K$ has finite index in $L$ and in $K$.
\end{proposition*}

{\it Proof}. We can assume that $\CA$ normalizes $K$ and $L$. Then the image of the group $K$ in the Lie group $\CA/L$ is a compact totally discontinuous subgroup in a Lie group and hence is finite.\qed
\bigskip 

We define the space of  test functions $\CS(\CA)$ on a moderate group $A$ as follows. The space $\CS(\CA)$
consists of complex valued functions $f$ on $\CA$ such that

\begin{enumerate}
	\item $f$ has compact support;

	\item
	$f$ is left invariant with respect to some basic compact subgroup;
	\item
	$f$ is a smooth function on the smooth manifold $K\setminus \CA$.
\end{enumerate}

   A function $f$ on $\CA$ is called smooth if in a neighborhood of
   any point it coincides with some test function. The algebra of 
   smooth functions will be denoted by $C^\infty(\CA)$.
   
 \subsubsection{Quotient spaces}  Let $X$ be a quotient space of $\CA$, i.e., $X= \CA/\CBB$ for a closed moderate subgroup $\CBB\subset \CA$ and  $X$ is endowed with the quotient topology.  We call such a space $X$ a moderate space. 
	
	We denote by $C^\infty(X)$ the algebra of functions that lift to smooth functions on $\CA$, and we denote
  by $\CS(X)$ the space of test functions on $X$, i.e., the space of smooth functions of compact support on $X$.

\begin{proposition*}Let $\al:\CA\to \CA'$ be a morphism of moderate groups, $X$, $X'$ quotients spaces of $\CA$ and $\CA'$, and $\bt:X\to X'$ a continuous map compatible with $\al$. Then $\bt$ is smooth, i.e., $\bt^*:C^\8(X')\to C^\8(X)$. 
\end{proposition*}

\subsubsection{Haar measure and co-invariants} Let $X$ be a moderate space. A Radon measure $\mu$ on $X$ defines a functional $I_\mu: \CS(X)\to \bc$, i.e., $I_\mu(f)=\int_X f d\mu$ for $f \in \CS(X)$.

For a torsor $L$, we can consider measures with values in $L$. Such a measure $\mu$ on $X$ defines a functional $I_\mu: \CS(X)\to L$.\\ 

{\bf Notation:} We denote by $L(\CA)=\CS_\CA(\CA):=\CS(\CA)/\langle f-a\circ f\rangle$ the space of co-invariants of the action of $A$ acting on the left   on $\CS(\CA)$ (i.e., $a\circ f(\al)= f(a\inv \al)$).

\begin{theorem*} Let $\CA$ be a moderate group. We consider the left action of $\CA$  on itself.

\begin{enumerate}
	\item The space of co-invariants $L=L(\CA)$ is a torsor. $\CA$ acts on $L$ trivially on the left and with some character $\Dl_\CA$ (the modulus character) on the right.  
	\item The canonical morphism $I: \CS(\CA)\to L$ is defined by a Radon measure $\mu_\CA$ with values in $L$. 
	\item The morphism $I$ and the measure $\mu_\CA$ are canonical, and it is  invariant with respect to  left and right actions of $\CA$ on $\CS(\CA)$. 
\end{enumerate}
\end{theorem*}

The theorem is essentially a reformulation of the Haar theorem. We call $\mu_\CA$ the Haar measure of $\CA$.  

\begin{remark} While the canonical map $I$  is defined initially only on the space of test functions, it could be extended to bigger spaces. Later we will apply $I$ also to some other classes of functions  using an appropriate regularization.
\end{remark}

We have the analogous construction for moderate quotient spaces. Let $X=\CA/\CBB$ be a quotient space of a moderate group $\CA$. Assume that there is a left $\CA$-invariant measure on $X$. The space $L(X)=\CS(X)_\CA$ of co-invariants is then a torsor, and there exists a canonical Haar measure $\mu_X$ on $X$ with values in $L(X)$ such that the map $I_{\mu_X}:\CS(X)\to L(X)$ is the canonical projection. 

\begin{proposition*}  We have the canonical isomorphism
$L(\CA)\simeq L(X)\otimes L(\CBB)$. 
\end{proposition*} The isomorphism is given by the integration (with values in co-invariants) along fibers. 
In particular, for a discrete subgroup $\CBB$, we have the canonical isomorphism $L(X)\simeq L(\CA)$, and hence the canonical integration map 
\begin{equation}\label{co-inv-iso}
	I_{X}:\CS(X)\to L(\CA)\ .
\end{equation}

\subsubsection{Groups over global fields}\label{glob-field-co-inv} Let $k$ be a global field. Let $\CG$ be an affine algebraic group defined over $k$. For every place $\fp$ of $k$, we consider the group of points $\CG_\fp=\CG(k_\fp)$ of $\CG$ over the local field $k_\fp$. We also consider the adelic group $\CG(\ba)$. It is defined with respect to compact open subgroups $\CG(\CO_\fp)\subset \CG(k_\fp)$ which are well-defined for almost all $\fp$. 

Let $\CV=\{V_\fp\}_{\fp\in\CP(k)}$ be a collection of complex vector spaces indexed by places of $k$. 

\begin{definition} An adelic structure $\Sigma$ on a family $\CV$ is a choice of vectors $v_\fp\in V_\fp$ for almost all $\fp$ (i.e., for all except finite number of places, up to a change of  vectors $v_\fp$ at finitely many places).
\end{definition}

\begin{definition} Let $\Sigma$ be an adelic structure on a family $\CV$. We define the restricted tensor product space $V$ by $V=\otimes_\Sigma V_\fp$.
\end{definition}
Namely, if $S\subset\CP(k)$ is a finite set then we define $V_S=\otimes_{\fp\in S}V_\fp$. If $S\subset S'$ and $S$ is sufficiently large, the adelic structure $\Sigma$ defines the canonical morphism $V_S\to V_{S'}$. Then, by the definition, $V=\otimes_\Sigma V_\fp=\lim\limits_{\stackrel{\rightarrow}{S}}V_S$.
\begin{remark} If all spaces $V_\fp$ are torsors and vectors $v_\fp$ are non-zero for almost all $\fp$, then $\otimes_\Sigma V_\fp$ is also  a torsor. 
\end{remark}


\begin{eexample*}
Let $\CG$ be an affine algebraic group defined over $k$.  For all $\fp$, we have the canonical map $I_\fp:\CS(\CG_\fp)\to L(\CG_\fp)$. 

\begin{claim} We have:

\begin{enumerate}
	\item The family of torsors $L(\CG_\fp)$ has canonical adelic structure $\Sigma_M$.
	\item There is  canonical isomorphism $L(\CG(\ba))\simeq \otimes_{\Sigma_M} L(\CG_\fp)$. 
\end{enumerate}
\end{claim}

Here the canonical adelic structure $\Sigma_M$ on $\{L(\CG_\fp)\}$ is obtained by taking the image $I_\fp(\chi_{K_\fp})$ of the characteristic function $\chi_{K_\fp}$ of the standard compact subgroup $K_\fp=\CG(\CO_\fp)$ at unramified places $\fp$. 
\end{eexample*}

\begin{remark} We note that in order to have the ``usual'' integral with respect to a measure with values in $\bc$, one has to choose isomorphisms $i_\fp: \{L(\CG_\fp)\}\simeq \bc$ for all places $\fp$, such that for almost all places, these satisfy $i_\fp(I_\fp(\chi_{K_\fp}))=1\in\bc$. This is easily translated into the familiar normalization of the local Haar measure by the standard compact subgroup. \end{remark}

\subsubsection{Tamagawa structure}\label{tam-meas-sect} There exists another remarkable adelic structure  $\Sigma_T$ for the family $\{L(\CG_\fp)\}$ proposed by T. Tamagawa \cite{Ta} (see also \cite{We}).

Let $A$ be an algebraic group defined over $k$.  We fix  a left invariant top differential form $\dl$ on $A$ defined over $k$. 
Such a choice gives rise to a measure $m(\dl_\fp)$ on  $A_\fp$, and in particular, defines the map $I_{m(\dl_\fp)}:\CS(A_\fp)\to\bc$ given by the integration. Hence we obtain the isomorphism $i_{m(\dl_\fp)}:L(A_\fp)\simeq \bc$ of the torsor of  co-invariants with the trivial torsor $\bc$. We can now define the Tamagawa adelic structure $\Sigma_T$ on the family  $\{L(\CG_\fp)\}$ by choosing the vector $t_\fp=i_{m(\dl_\fp)}\inv(1)\in L(\CG_\fp)$ for all $\fp$. We call the resulting restricted tensor product torsor $L^T(\CG(\ba))\simeq \otimes_{\Sigma_T} L(\CG_\fp)$ the Tamagawa torsor. Note that since non-zero vectors $t_\fp$ are specified for {\it all} places $\fp$, the torsor $L^T(\CG(\ba))$ comes with the canonical  trivialization given by the ``Tamagawa measure'', i.e., by the vector  $\ft=\ft_\dl=\otimes_\fp t_\fp$. The Tamagawa structure $\ft_\dl$ does not depend on the rational class of the form $\dl$ as follows from the standard product formula.

 \begin{remark}We do not claim that  torsors $L^T(\CG(\ba))$ and $L(\CG(\ba))$ are isomorphic with respect to a collection of some local isomorphisms $j_\fp:L(\CG_\fp)\to L(\CG_\fp)$ mapping the adelic structure $\Sigma_M$ to $\Sigma_T$ at almost all places. If this is the case, one can integrate functions in $\CS(\CG(\ba))$ with respect to the Tamagawa structure $\ft$ (i.e., this gives the usual Tamagawa measure). 
Sometimes such an isomorphism exists and it is possible to integrate functions in $\CS(\CG(\ba))$ with respect to $\ft$ (e.g., for a unipotent subgroup $N\simeq k^+$), and this means that the Tamagawa construction provides  a measure in the usual sense.  However, in general, we cannot integrate functions in $\CS(\CG(\ba))$ with respect to $\ft$ since the Euler product $\prod_{\fp\not\in S}
i_{|\dl_\fp|}(I_\fp(\chi_{K_\fp}))$ is not absolutely convergent (e.g., for the  torus $A\simeq k^\times$). This appears when two local trivializations $(L(\CG_\fp),I_\fp(\chi_{K_\fp}))\simeq (\bc, 1)$ and $(L(\CG_\fp),t_\fp)\simeq (\bc, 1)$ are not globally compatible, and one has to introduce a regularization procedure in order to obtain a measure out of the Tamagawa structure $\ft$ (i.e., to construct another trivialization of $L^T(\CG(\ba))$). \end{remark}

\subsubsection{Characters}\label{char} We also consider integration  twisted by characters. 

Let $(V_\tau, \tau)$ be a representation of $\CA$, and $\chi:\CA\to\bc^\times$ be a character. We have the Jacquet module $J_\chi(\pi)=V_\tau/\langle v-\chi\inv(a)\tau(a) v\rangle$, $v\in V_\tau$.

Let $X= \CA/\CBB$ be a homogenous $\CA$-space.   We denote by $L_\chi(X)=J_\chi(\CS(X))$ the corresponding Jacquet module. Let us assume that on $X$ there is an invariant measure. We can  describe this torsor as follows. Let $C(X,\chi)$ be the space of functions on $X$ satisfying $f(ax)=\chi(a)f(x)$. This space is zero  if $\chi|_\CBB\not\equiv 1$, and is a torsor otherwise.  
\begin{claim}There is a canonical isomorphism $L_\chi(X)\simeq C(X,\chi)\otimes L(X)$.
\end{claim}
Choice of a point $x_0\in X$ gives a trivialization $C(X,\chi)\simeq \bc_\chi$, and hence the isomorphism $L_\chi(X)\simeq L(X)\otimes\bc_\chi$.
Hence $L_\chi(\CA)$ is a torsor on which $\CA$ acts by $\chi$ on the left and by $\Dl_\CA\chi\inv$ on the right.

The natural projection $I_\chi:\CS(X)\to L_\chi(X)$ corresponds to  the integration with some measure $\mu_{(X,\chi)}$ with values in $L_\chi(X)$. 



Let $\CG$ be an affine algebraic group  defined over $k$. Let $\chi$ be a character $\chi=\otimes_\fp\chi_\fp$ of $\CG(\ba)$.
\begin{claim} We have the isomorphism $L_\chi(\CG(\ba))\simeq \otimes_{\Sigma_M}L_{\chi_\fp}(\CG_\fp)$.
\end{claim}

Consider the automorphic space $X_\CG=\CG(k)\setminus \CG(\ba)$.  Let $\chi$ be a character of $\CG(\ba)$ which is trivial on $\CG(k)$ (i.e., $\chi: \CG(k)\setminus \CG(\ba)\to\bc$).
We can trivialize the torsor $C(X_\CG,\chi)$ by using the evaluation at the base point $x_0=\{\CG(k)\}\in X_\CG$. This gives the isomorphism 
\begin{equation} L_\chi(X_\CG)\simeq L(X_\CG)\otimes\bc_\chi\simeq L_\chi(\CG(\ba))\simeq \otimes_{\Sigma_M}L_{\chi_\fp}(\CG_\fp)\ .
\end{equation}

As a result,  we  have the corresponding integration map $I_{X_\CG,\chi}:\CS(X_\CG)\to L(\CG(\ba))\otimes\bc_\chi$.

\subsection{Automorphic representations} Let $\CG$ be a reductive algebraic group defined over $k$. Let $\pi$ be an irreducible smooth  representation of the adelic group $\CG(\ba)$. We denote by $V_\pi$ the space of smooth vectors of $\pi$ and by $\om_\pi$ the central character of $\pi$. We have decompositions  $\pi=\otimes_\fp\pi_\fp$ and $V_\pi=\hat\otimes_\fp V_{\pi_\fp}$ into the restricted tensor product of local representations.  

Let $X_{\CG}=\CG(k)\setminus \CG(\ba)$ be the automorphic space. An automorphic structure on an (abstract) adelic representation $\pi$ is an intertwining  map  $\nu:V_\pi\to \CF(X_{\CG})$  with the representation of $\CG(\ba)$ in the space   of  functions on $ X_{\CG}$.  We call a pair $(\pi,\nu)$ an automorphic  representation.
For a cuspidal $(\pi,\nu)$, the image of $\nu$ belongs to the space of rapidly decreasing smooth functions on $X_{\CG}$. For a vector $v\in V_\pi$, we will denote by $\phi_v=\nu(v)\in \CG(\ba)$ the  corresponding automorphic function.


     We denote by $S(\pi)$ the set of places (including infinite places) where  $\pi$ is ramified (i.e., the complement to the set of unramified places $\fp$ where the standard $K_\fp$-fixed vector $e^0_\fp\in V_{\pi_\fp}$ is specified).

\subsection{Periods}\label{periods-setup} Let $\CH\subset \CG$ be an algebraic subgroup defined over $k$. Denote by  $ X_{\CH}=\CH(k)\setminus \CH(\ba)\subset X_{\CG}$ the closed $\CH(\ba)$-orbit.  Let $\chi: \CH(k)\setminus \CH(\ba)\to \bc$ be a character of $\CH(\ba)$ which is trivial on $\CH(k)$.   According to \eqref{co-inv-iso}, we have the integration map  $I_{X_{\CH},\chi}:\CS(X_{\CH})\to L_\chi(\CH(\ba))$. This together with the automorphic realization map $\nu$ and the restriction map $res_{X_{\CH}}:C^\8(X_G)\to C^\8(X_{\CH})$ give rise to the $\CH(\ba)$-equivariant period map $p_{\CH,\chi}=I_{X_{\CH},\chi}\circ res_{\fO_{\CH}}\circ\nu:
V_\pi\to L_\chi(\CH(\ba))\simeq\otimes_{\Sigma_M} L_\chi(\CH_\fp)$, where $\CH_\fp=\CH(k_\fp)$.  Formally, we need to assume that the corresponding integrals are well-defined (e.g., the orbit $X_{\CH}$  is compact or the automorphic representation $(\pi,\nu)$ is cuspidal).

\begin{definition}
\hfill
\begin{enumerate}
	\item The space  $P(V_\pi,L_\chi(\CH(\ba))):=\Hom_{H(\ba)}(V_\pi, L_\chi(\CH(\ba)))$  is called the period space. 
	\item For every place $\fp$, the space $P(V_{\pi_\fp},L_\chi(\CH_\fp)):=\Hom_{\CH_\fp}(V_{\pi_\fp}, L_\chi(\CH_\fp))$ is called the local period space.
\end{enumerate}
\end{definition}

With the decomposition $\chi=\otimes\chiv$, we have the factorization  \begin{equation*}P(V_\pi,L_\chi(\CH(\ba)))\simeq \hat\otimes_\fp P(V_{\pi_\fp},L_\chiv(\CH_\fp))\ .
\end{equation*}

We  will  assume that the local period space $P(V_{\pi_\fp},L_\chiv(\CH_\fp))$ is at most one-dimensional. 
Hence any map in the period space  is  factorizable, and we can choose a factorization for the automorphic period $p_{\CH,\chi}$. To choose a factorization of the torsor $P(V_\pi,L_\chi(\CH(\ba)))$ into a restricted tensor product, we need to choose for almost all  places $\fp$, a special vector $p_\fp^0\in P(V_{\pi_\fp},L_\chiv(\CH_\fp))$. We choose it  by requiring that $ p^0_\fp(e^0_\fp)=I_{\CH_\fp,\chiv}(\chi_{K_{\CH_\fp}})$ (in fact, one has to check that such a normalization is possible, i.e., that there exists a non-vanishing invariant map on the standard vector $e^0_\fp$). 
 Then for a sufficiently large finite set $S\subset\CP(k)$, we have $p_{\CH,\chi}=(\otimes_{\fp\in S}p_\fp)\otimes(\otimes_{\fp\not\in S}p^0_\fp)$ for some choice of local ramified components $p_\fp$ for $\fp\in S$.

\subsection{Action on periods}\label{idea-coinv} We reformulate our scheme from Section \ref{idea} in the language of co-invariants.
Let $\CH_1, \CH_2\subset\CG$ be two algebraic subgroups as above. In particular, we will assume that all local spaces satisfy the Gelfand condition of multiplicity one. 

\subsubsection{Local maps}\label{loc-map-coinv} Let $\fp$ be a place of $k$, $(\pi_\fp, V_\fp)$ be an irreducible  smooth  representation of $\CG_\fp=\CG(k_\fp)$. Let $\chi_{i,\fp}:\CH_{i,\fp}\to\bc$ be  a character of $\CH_{i,\fp}$. We fix a non-zero  invariant differential form $\dl$ on $\CH_1$ defined over $k$. Let $\dl_{\fp}$ be the corresponding invariant  measure on $\CH_{1,\fp}$.  We use the measure $\dl_\fp$ to trivialize by the map $\ev_{\dl_\fp}\!:
L_{\chi_{1,\fp}}(\CH_{1,\fp})\MapsTo \bc_{\chi_{1,\fp}}$  the corresponding space of co-invariants, and, correspondingly, we get the isomorphism (see Section \ref{char})
\begin{equation*}
	\ev_{\dl_\fp}^*\!: P(V_{\pi_\fp},L_{\chi_{1,\fp}}(\CH_{1,\fp}))\MapsTo \Hom_{A_\fp}(V_{\pi_\fp},\bc_{\chi_{1,\fp}})\ .
\end{equation*}   
We also consider the integration map 
 $I_{\CH_{2,\fp},\chi_{2,\fp}}:\CS(\CH_{2,\fp})\to L_{\chi_{2,\fp}}(\CH_{2,\fp})$. 

We now construct the local map between local period spaces
\begin{align}\label{i-psi-chi-local-intro}
\fii(\chi_{1,\fp}, \chi_{2,\fp}, \dl_\fp): P(V_{\pi_\fp},L_{\chi_{1,\fp}}(\CH_{1,\fp}))\to P(V_{\pi_\fp},L_{\chi_{2,\fp}}(\CH_{2,\fp}))
\end{align} 
using  the co-invariant map $I_{\CH_{2,\fp}}$ and the evaluation map $\ev^*_{\dl_\fp}$. Namely, for a given local period  $p_{(\CH_{1,\fp},\chi_{1,\fp})}\in P(V_{\pi_\fp},L_{\chi_{1,\fp}}(\CH_{1,\fp}))$, we consider the function $f:\CH_{2,\fp}\to \Hom_{\CH_{1,\fp}}(V_{\pi_\fp},\bc_\psiv)\subset V^*_\piv$ given by $f(h_2)=\piv^*(h_2)(\ev^*_{\dl_\fp}(p_{(\CH_{1,\fp},\chi_{1,\fp})}))$  (here $\piv^*$ is the dual to $\piv$ representation). We now consider the matrix coefficient function  $f_v(h_2)=\piv^*(h_2)(\ev^*_{\dl_\fp}(p_{(\CH_{1,\fp},\chi_{1,\fp})}))(v)$, $v\in V$, and apply to the function $f_v(h_2)$ on $\CH_{2,\fp}$ the integration map $I_{\CH_{2,\fp},\chi_{2,\fp}}(f_v)\in L_{\chi_{2,\fp}}(\CH_{2,\fp})$.  Formally, the function $f_v$ is not compactly supported and we have to make sense of the corresponding integral. This is achieved by considering appropriate regularization procedure (e.g., by the analytic continuation method).

\subsubsection{Global maps} In order to define the global map 
\begin{align}\label{i-psi-chi-global-intro}
\fii(\chi_{1}, \chi_{2}, \dl): P(V_{\pi},L_{\chi_{1}}(\CH_{1}))\to P(V_{\pi},L_{\chi_{2}}(\CH_{2}))\ ,
\end{align} 
we  make sense out of the Euler product $\hat{\otimes}_\fp \fii(\chi_{1,\fp}, \chi_{2,\fp}, \dl_\fp)$ using methods we described in Section \ref{idea}.

\section{Whittaker and Hecke periods relations}\label{Hecke-sect} We want to illustrate how the  procedure described  in Section \ref{idea-coinv}  relates Whittaker and Hecke functionals on an automorphic cuspidal representation of $G=\GL(2)$.

\subsection{Whittaker and Hecke periods}\label{Wh-Hecke}
\subsubsection{Whittaker period} We fix a nontrivial additive character $\psi:\ba\to \bc^\times$ which is trivial on principal adeles. We  view the character $\psi$ as a character of the unipotent group $N(\ba)$. We consider the $N(\ba)$-orbit $X_N=N(k)\setminus N(\ba)\subset X_G$ and the corresponding period it induces on an automorphic  representation $(\pi,\nu)$ of $G=GL(2)$. 

According to the above scheme, we  view the Whittaker period $p_{(N,\psi)}$ as an element in the period space $P(V_\pi,L_\psi(N(\ba)))=\Hom_{N(\ba)}(V_\pi,L_\psi(N(\ba)))$. In the factorization of the Whittaker period $p_{(N,\psi)}=\otimes_\fp p_{(N_\fp,\psi_\fp)}$  for almost all $\fp$, the local component $p_{(N_\fp,\psi_\fp)}\in P(V_{\pi_\fp},L_{\psi_\fp}(N_\fp))= \Hom_{N_\fp}(V_{\pi_\fp},L_{\psi_\fp}(N_\fp))$ is unramified (i.e., maps the standard vector $e^0_\fp$ to the image in co-invariants of the characteristic function  $\chi_{N(\CO_\fp)}$ of the set $N(\CO_\fp)$).

\begin{remark} In more classical terms, the Whittaker period/functional on $V_\pi$ is given by the integral
\begin{align}\label{whittak-def}
\CW^\psi(v)=\int_{N(k)\setminus N(\ba)}\psi\inv (n)\phi_v(n)\ dn\ .	
\end{align}
Here $dn$ is the  measure on $N(\ba)$ obtained from an invariant differential form and $\phi_v$ is the automorphic function corresponding to the vector $v\in V_\pi$. We have $\CW^\psi\in \Hom_{N(\ba)}(V_\pi,\bc_\psi)$, where $\bc_\psi$ is the one-dimensional  $N(\ba)$-module  with the action given by $\psi$. It is well-known that $\dim \Hom_{N(\ba)}(V_\pi,\bc_\psi)=1$,  the space of local functionals $\Hom_{N_\fp}(V_{\pi_\fp},\bc_{\psi_\fp})$ is also one-dimensional, and the global space decomposes into the restricted product of local spaces. We have the following standard decomposition of the automorphic Whittaker functional $\CW^{\psi}$  into a product of local functionals. For an unramified place $\fp\not\in S(\pi,\psi)$ (here $S(\pi,\psi)$ denotes the set of primes where $\pi$ or $\psi$ are ramified), let $\CW_0^{\psi_\fp}\in \Hom_{N_\fp}(V_{\pi_\fp},\bc_{\psi_\fp})$ be the local functional satisfying $\CW_0^{\psi_\fp}(e^0_\fp)=1$ for the standard $K_\fp$-fixed vector $e^0_\fp\in V_{\pi_\fp}$. We then choose local functionals $\widetilde\CW^{\psi_\fp}\in \Hom_{N_\fp}(V_{\pi_\fp},\bc_{\psi_\fp})$ for ramified primes, so that $\CW^{\psi}=\otimes_{\fp\in S(\pi, \psi)}\widetilde\CW^{\psi_\fp}\otimes_{\fp\not\in S(\pi, \psi)}\CW_0^{\psi_\fp}$.
\end{remark}

\subsubsection{Hecke period} We now consider the Hecke period. Let $\chi$ be a Hecke (quasi-)character of $k$ and let $\chi: A(\ba)\to \bc^\times$,  $\chi\left( \left(\begin{smallmatrix}
a & \\
                         &    1      \\
\end{smallmatrix} \right)\right)=\chi(a)$,
be the corresponding  (quasi-)character of $A(\ba)$ trivial on the principal subgroup $A(k)$. We consider the orbit $X_A=A(k)\setminus A(\ba)\subset X_G$, and the corresponding (Hecke) period  $d_\chi$ it induces on a cuspidal automorphic representation $(\pi,\nu)$.
According to the above scheme, we can  view the Hecke period as an element
in the torsor of periods $P(V_\pi,L_\chi( A(\ba)))=\Hom_{A(\ba)}(V_\pi, L_\chi(A(\ba)))$. We  have a factorization
 $p_{(A,\chi)}=\otimes_\fp p_{(A_\fp,\chi_\fp)}$ where for almost all $\fp$, the local component $p_{(A_\fp,\chi_\fp)}\in P(V_{\pi_\fp},L_{\chi_\fp}(A_\fp))$ is unramified (i.e., maps the standard vector $e^0_\fp$ to the image in co-invariants of the characteristic function  $\chi_{A(\CO_\fp)}$ of the set $A(\CO_\fp)$).


\begin{remark} In more classical terms, we have the following description of the functional   $d_\chi=d_\chi(\pi): V_\pi\to\bc_\chi$ (here $\bc_\chi$ denotes the one-dimensional $A(\ba)$-module with the action  given by $\chi$). We fix an invariant rational differential form on $A$ and denote by $\mu$ the corresponding invariant measure. As an element in the space $\Hom_{A(\ba)}(V_\pi, \bc_\chi)$, the corresponding period functional is given by the integral
\begin{align}\label{torus-per-def}
d_\chi(v)=\int_{A(k)\stm A(\ba)}\chi\inv(\bar a)\phi_v(\bar a)\ \mu\ ,
\end{align}
for $v\in V_\pi$ and $\phi_v$ the corresponding  automorphic function. The integral is absolutely convergent since functions in a cuspidal representation are rapidly decreasing at infinity. It is well-known that $\dim  \Hom_{A(\ba)}(V_\pi,\bc_\chi)=1$, that the local space $\Hom_{A_\fp}(V_{\pi_\fp},\bc_{\chi_\fp})$ also satisfies the multiplicity one property, and hence $ \Hom_{A(\ba)}(V_\pi,\bc_\chi)=\hat\otimes_\fp  \Hom_{A_\fp}(V_{\pi_\fp},\bc_{\chi_\fp})$. We again can choose a decomposition  $d_\chi=\otimes_{\fp\in S(\pi, \chi)}\widetilde d_{\chi_\fp}\otimes_{\fp\not\in S(\pi, \chi)}d^0_{\chi_\fp}$ into local components with $d^0_{\chi_\fp}(e^0_\fp)=1$. In what follows we recall Hecke-Jacquet-Langlands recipe
how to choose local components consistently at all places.  

\end{remark}

\subsection{Action from Whittaker to Hecke}\label{W-H-sec}

Let us describe a morphism from the Whittaker torsor to the Hecke torsor.

\subsubsection{Local map} Let $\fp$ be a place of $k$ and $(\pi_\fp, V_\fp)$ be an irreducible  smooth  representation of $G_\fp=G(k_\fp)$. Let $\psi_\fp:N_\fp\to\bc$ be a nontrivial  character of $N_\fp$ and $\chi_\fp:A_\fp\to\bc$ be a character of $A_\fp \simeq k_\fp^\times$. We fix a (non-zero) invariant differential form $\dl_N$ on $N$ and consider the corresponding invariant measure $dn_\fp=dn_\fp(\dl_N)$   on $N_\fp$.  The measure $dn_\fp$ gives rise to the trivialization  $\ev_{dn_\fp}:
L_{\psi_\fp}(N_\fp)\MapsTo \bc_\psiv$ of the corresponding co-invariants and to  the isomorphism $\ev_{dn_\fp}^*: P(V_{\pi_\fp},L_{\psi_\fp}(N_\fp))\MapsTo \Hom_{N_\fp}(V_{\pi_\fp},\bc_\psiv)$. We  consider the integration map $I_{A_\fp,\chi_\fp}:\CS(A_\fp)\to L_\chiv(A_\fp)$ (see Section \ref{char}). 

Following the scheme formulated in Section \ref{idea-coinv}, we consider  the local map 
\begin{align}\label{i-chi-psi-local}
\fii(\chi_\fp,\psi_\fp, dn_\fp): P(V_{\pi_\fp},L_{\psi_\fp}(N_\fp))\to P(V_{\pi_\fp},L_{\chi_\fp}( A_\fp))
\end{align} 
constructed out of the co-invariant map $I_{A_\fp,\chi_\fp}$ and out of the evaluation map $\ev^*_{dn_\fp}$. The map $\fii(\chi_\fp,\psi_\fp, dn_\fp)$ could be described in the following terms. For a map $p_{(N_\fp,\psi_\fp)}\in P(V_{\pi_\fp},L_{\psi_\fp}(N_\fp))$ and a vector $v\in V_\fp$, we consider the  (matrix coefficient) function $\al_v(\bar a_\fp)= \ev_{dn_\fp}(p_{(N_\fp,\psi_\fp)}(\pi_\fp(\bar a_\fp)v))\in C^\8(A)$ and then take its image $I_{A_\fp,\chi_\fp}(\al_v)\in L_{\chi_\fp}(A_\fp)$ under the integration, i.e., 
\begin{align}\label{p-T-def} p_{(A_\fp,\chi_\fp)}(v)=[\fii(\chi_\fp,\psi_\fp,dn_\fp)(p_{(N_\fp,\psi_\fp)})](v):=   I_{A_\fp,\chi_\fp}(\al_v)\ . \end{align}

Smooth functions $\al_v(\bar a_\fp)$ with $v\in V_\fp$ obtained in such a way are not compactly supported on $A_\fp$. Hence, in fact, we have to extend the integration map $I_{A_\fp,\chi_\fp}:\CS(A_\fp)\to L_{\chi_\fp}(A_\fp)$ to such functions (i.e., to the space of matrix coefficients $\al_v$ as above for an irreducible representation $\pi_\fp$). 


\begin{proposition*}\label{loc-i-map-prop}\ 
\begin{enumerate}
	\item For $Re(\chi_\fp)\ll 1$, the map $\fii(\chi_\fp,\psi_\fp, dn_\fp)$ is well-defined (by an absolutely convergent integral \eqref{i-chi-psi-local-int}). It has the  meromorphic continuation to the complex space of all characters (i.e., to the complex plane of characters of the form $\chi_\fp|\cdot|^{-s}_\fp$).
	\item For the unramified  data, we obtain the Hecke-Jacquet-Langlands local $L$-factor. Namely, 
	\begin{align}\label{L-fac-co-inv}\fii(\chi_\fp,\psi_\fp, dn_\fp)(\CW_0^{\psi_\fp})=L_\fp(\chiv,\pi_\fp)\cdot d_\chiv^0\  .
	\end{align}

\end{enumerate}
\end{proposition*}
The meaning of the unramified condition above is as follows (see Section \ref{Wh-Hecke}): $\pi_\fp$ is an unramified representation, $\psi_\fp$ has conductor  $cond(\psiv)=\CO_\fp$, $\chi_\fp$ is an unramified character, and the form $\dl_N$ comes from a rational invariant form on $N$ which is unramified at $\fp$ (and as a result satisfies  $ dn_\fp(\dl_N)(N(\CO_\fp))=1$), the unramified Whittaker functional $\CW_0^{\psi_\fp}\in P(V_\piv,L_\psiv(N_\fp))$ satisfies $\CW_0^{\psi_\fp}(e^0_\fp)=l^0_\psiv(N_\fp)$ for $l^0_\psiv(N_\fp)\in L_\psiv(N_\fp)$ given by the adelic structure on  $L_\psiv(N_\fp)$ described in Section \ref{glob-field-co-inv}, and correspondingly for the  functional $d_0^{\chi_\fp}\in P(V_\piv,L_\chiv(A_\fp))$ with $d_0^{\chi_\fp}(e^0_\fp)=l^0_\chiv(A_\fp)$ and $l^0_\chiv(A_\fp)\in L_\chiv(A_\fp)$. 
\begin{remark} 

We claim that the map $\fii(\chi_\fp,\psi_\fp, dn_\fp)$ in \eqref{i-chi-psi-local}
naturally appears, in another language, in \cite{JL} as local zeta integrals on $\GLt$ of Jacquet and Langlands.  
Let us fix a non-zero invariant local measure $d^\times a_\fp$ on $A_\fp$. This gives rise to isomorphisms $\ev_{d^\times a_\fp}:L_\chi(A_\fp)\MapsTo \bc_\chiv$ and  $\ev_{d^\times a_\fp}^*:P(V_{\pi_\fp},L_{\chi_\fp}(A_\fp))\MapsTo  \Hom_{A_\fp}(V_\pi,\bc_{\chi_\fp}) $.
Using isomorphisms  $\ev_{dn_\fp}^*$ and $\ev_{d^\times a_\fp}^*$,  we see that the map \eqref{i-chi-psi-local}   could be described by the following standard in $\GL(2)$ theory integral:
\begin{align}\label{i-chi-psi-local-int}
[i(\chi_\fp,\psi_\fp, dn_\fp, d^\times a_\fp)(p_{(N_\fp,\psi_\fp)})](v):&=&\\ [\ev_{d^\times a_\fp}^*( p_{(A_\fp,\chi_\fp)})](v)&=&\int_{A_\fp}\chi_\fp\inv(\bar a)l^{\psi_\fp}(\pi_\fp(\bar a)v) d^\times a_\fp \ ,\nonumber
\end{align} 
for $v\in V_{\pi_\fp}$ and  $l^{\psi_\fp}= \ev_{dn_\fp}^*(p_{(N_\fp,\psi_\fp)})\in \Hom_{N_\fp}(V_\piv,\bc_\psiv)$. Translated into local zeta integrals \eqref{i-chi-psi-local-int}, the relation \eqref{L-fac-co-inv} reads 
\begin{align*}[i(\chiv,\psi_\fp,dn_\fp,d_1^\times a_\fp)(\CW_0^{\psi_\fp})](e^0_\fp)=L_\fp(\chiv,\pi_\fp)\ ,
	\end{align*} for the measure $d_1^\times a_\fp$ on $A_\fp$ normalized  by the condition  $d_1^\times a_\fp(A(\CO_\fp))=1$.

Note that while the integral \eqref{i-chi-psi-local-int} depends on the choice of the measure $d^\times a_\fp$, the map \eqref{i-chi-psi-local} does not. Over archimedian fields, a  similar approach  appeared in \cite{Po}. 


\end{remark}

\subsubsection{Global map} Fix an automorphic cuspidal representation $(\pi,\nu)$ and a non-trivial  character $\psi: N(k)\stm N(\ba)\to\bc$. Choose an invariant  differential form $\dl_N$ on $N$.  We want to define a map $\fii(\chi,\psi, \dl_N): P(V_{\pi},L_\psi(N(\ba)))\to P(V_{\pi},L_\chi(A(\ba)))$ as a tensor product of local maps. 
\begin{proposition*}\label{glob-i-map-prop}  Fix an invariant rational differential form $\dl_{N}$ on $N$.

\begin{enumerate}
	\item  The tensor product $\fii(\chi,\psi, \dl_N)=\otimes_\fp \fii(\chi_\fp,\psi_\fp, dn_\fp)$ is absolutely convergent for $Re(\chi)\ll 1$, and has the meromorphic continuation to the complex space of all characters. 
	\item The resulting map 
\begin{align*}\label{i-chi-psi-global}
\fii(\chi,\psi)=\fii(\chi,\psi, \dl_N): P(V_{\pi},L_\psi(N(\ba)))\to P(V_{\pi},L_\chi(A(\ba)))\ 
\end{align*}
does not depend on the choice of the  form $\dl_N$. 
\end{enumerate}
\end{proposition*}


 \subsubsection{Action on automorphic periods} We now came to the last step of our scheme where we compute the effect of the defined map on automorphic periods. 
\begin{theorem*}\label{N-per-T-per-thm}  The global map $\fii(\chi,\psi)$ is coherent, i.e., it sends the automorphic Whittaker  period $\CW^\psi$ to the automorphic Hecke period $d_{\chi}$. Namely, we have \[ \fii(\chi,\psi)(\CW^\psi)=d_{\chi}\ .\]

\end{theorem*}

\begin{remarks} \label{rem-bundls-hol-sect}
1. The set of Hecke quasi-characters $\chi: A(k)\sm A(\ba)\to \bc^\times$ has the natural structure of a complex manifold (with infinitely many connected components).  For a given quasi-character $\chi$ we denote by $C_\chi:=\{\chi_s=\chi|\cdot|^s,\ s\in\bc\}\simeq \bc$ the corresponding leaf. The set $\widehat{N(k)\sm N(\ba)}$ is discrete and we denote by $pt_\psi$ the point corresponding to a character $\psi$. For  given character $\psi : N(k)\sm N(\ba)\to\bc$ and a quasi-character $\chi: A(k)\sm A(\ba)\to \bc^\times$, consider the trivial line bundle $W_{\psi,\chi}(\pi)$ over $pt_\psi \times C_\chi\simeq C_\chi$,  with the constant fiber $P_{\pi,\psi}=P(V_{\pi},L_\psi(N(\ba)))$, and the line bundle $H_\chi(\pi)$ over $C_\chi$ with the fiber $P(V_{\pi},L_{\chi_s}(A(\ba)))$. Both line bundles have natural  holomorphic (for a cuspidal representation $\pi$) sections $S^{aut}_\psi=\CW^\psi$ and $S^{aut}_{\chi_s}=d_\chi$ coming from automorphic periods (the section for the Whittaker bundle is a constant non-zero section providing a trivialization, but the automorphic Hecke section of $H_\chi(\pi)$ has zeroes). 
 
One can naturally view the map $\fii(\chi,\psi)$ as a holomorphic map between two bundles which is trivial on the base and which maps the section  $S^{aut}_\psi$ to the section $S^{aut}_{\chi_s}$.

2. Our treatment of the map $\fii(\chi,\psi)$ is based on the possibility to include it in an analytic family of maps $\fii(\chi_s,\psi)$ and use analytic continuation. One can use the single character $\chi$ and follow the recipe formulated in the second part of Section \ref{idea} which is applicable when there are no deformations. Namely, by looking at the local unramified computation \eqref{L-fac-co-inv}, we see that the partial Euler product $\otimes_{\fp\not\in S} L_\fp(\chiv,\pi_\fp)\inv \fii(\chi_\fp,\psi_\fp, dn_\fp)$ is absolutely convergent (since  all but finite number of terms are trivial). Assuming that $L_S(\chi|\cdot|^s,\pi)$ have no pole at $s=0$ (e.g., $\pi$ is cuspidal), the formula \eqref{L-reg-intro} defines the map $\fii(\chi,\psi)(v)= L_S(\chi,\pi)\prod_{\fp\not\in S} d_\fp^0(v_\fp)\prod_{\fp\in S} d_\fp(v_\fp)$. Note that this map might be zero. Hence we can define the map $\fii(\chi,\psi)$ without considering the whole family $\fii(\chi_s,\psi)$ (of course this is a purely cosmetic change since we use the equivalent property of $L_S(\chi|\cdot|^s,\pi)$). What we can not see that way is the fact that the resulting map is coherent (i.e., Theorem \ref{N-per-T-per-thm}), and in fact we will see in Section \ref{conv-Hecke-sect} that for other periods this does not always holds. To show that $\fii(\chi_s,\psi)$ is coherent, we use the Hecke-Jacquet-Langlands unfolding which does not have an analog for other periods (e.g., for the map from Hecke period to Whittaker period considered in the next section).

\end{remarks}

\subsection{Action from Hecke to Whittaker}\label{conv-Hecke-sect}

Here we treat the opposite direction which inspite of similarities is not classical.

\subsubsection{Local map} Let $\fp$ be a place of $k$ and $(\pi_\fp, V_\fp)$ be an irreducible  smooth  representation of $G_\fp=G(k_\fp)$. Let $\psi_\fp:N_\fp\to\bc$ be a nontrivial  character of $N_\fp$ and $\chi_\fp:A_\fp\to\bc$ be a character of $A_\fp$. We fix a non-zero  invariant differential form $\dl_{A}$ on $A$ defined over $k$. Let $d^\times a_\fp=d^\times a_\fp(\dl_{A})$ be the corresponding invariant  measure on $A_\fp$.  We use the measure $d^\times a_\fp$ to trivialize $\ev_{d^\times a_\fp}:
L_{\chiv}(A_\fp)\MapsTo \bc_\chiv$  the corresponding co-invariants, and, correspondingly, we get the isomorphism $\ev_{d^\times a_\fp}^*: P(V_{\pi_\fp},L_{\chiv}(A_\fp))\MapsTo \Hom_{A_\fp}(V_{\pi_\fp},\bc_\chiv)$ (see Section \ref{char}). We also consider the integration map 
 $I_{N_\fp,\psiv}:\CS(N_\psiv)\to L_{\psiv}(N_\fp)$. 

We now construct the local map 
\begin{align}\label{i-psi-chi-local}
\fii(\psi_\fp, \chiv, d^\times a_\fp): P(V_{\pi_\fp},L_{\chiv}(A_\fp))\to P(V_{\pi_\fp},L_{\psiv}(N_\fp))
\end{align} 
using   maps $I_{N_\fp,\psiv}$ and $\ev^*_{d^\times a_\fp}$ as in the previous case. Namely, for a map $p_{(A_\fp,\chiv)}\in P(V_{\pi_\fp},L_{\chiv}(A_\fp))$ and a vector $v\in V_\fp$, we consider the  (matrix coefficient) function $\bt_v(n_\fp)= \ev_{d^\times a_\fp}(p_{(A_\fp,\chiv)}(\pi_\fp(n_\fp)v))\in C^\8(N_\fp)$ and then take its image $I_{N_\fp,\psiv}(\bt_v)\in L_{\psiv}(A_\fp)$ under the integration map, i.e., 
\begin{align}\label{p-N-def} p_{(N_\fp,\psiv)}(v)=[\fii(\psi_\fp,\chi_\fp,d^\times a_\fp)(p_{(A_\fp,\chiv)})](v):=   I_{N_\fp,\psiv}(\bt_v)\ . \end{align}

Smooth functions $\bt_v(n_\fp)$ with $v\in V_\fp$ obtained in such a way are not compactly supported on $N_\fp$. Hence we have to extend the integration map $I_{N_\fp,\psiv}:\CS(N_\fp)\to L_{\psiv}(N_\fp)$ to such functions (i.e., to the space of matrix coefficients $\bt_v$ as above for an irreducible representation $\pi_\fp$). Indeed it is easy to see that such an extension exists and is unique as explained in Section \ref{reg-whitak-f}.

\begin{proposition*}\label{loc-i-N-map-prop}


	 For the unramified  data, we have  
	\begin{align}\label{L-fac-Heck-converse}\fii(\psi_\fp,\chi_\fp, d^\times a_\fp)( d_\chiv^0)=L_\fp(\chiv,\pi_\fp)\inv\cdot \CW_0^{\psi_\fp}\ 
	\end{align}
with the Hecke-Jacquet-Langlands local $L$-factor $L_\fp(\chiv,\pi_\fp)$.

\end{proposition*}

For the meaning of the unramified data $d_\chiv^0$ and $\CW_0^{\psi_\fp}$, see Section \ref{Wh-Hecke}.

\begin{remark} 
Let us fix a non-zero invariant local measure $dn_\fp$ on $N_\fp$. This gives rise to isomorphisms $\ev_{dn_\fp}:L_\psiv(N_\fp)\MapsTo \bc_\psiv$ and  $\ev_{dn_\fp}^*:P(V_{\pi_\fp},L_{\psiv}(N_\fp))\MapsTo  \Hom_{N_\fp}(V_\pi,\bc_{\psiv}) $.
Using isomorphisms  $\ev_{dn_\fp}^*$ and $\ev_{d^\times a_\fp}^*$,  we see that the map \eqref{i-psi-chi-local}   could be described by the following  integral:
\begin{align}\label{i-psi-chi-local-int}
[i(\psi_\fp,\chi_\fp, d^\times a_\fp, dn_\fp)(p_{(A_\fp,\chiv)})](v):&=&\\ [\ev_{dn_\fp}^*( p_{(N_\fp,\psiv)})](v)&=&\int_{N_\fp}\psi_\fp\inv(n)d_{\chi_\fp}(\pi_\fp(n_\fp)v)\ dn_\fp \ ,\nonumber
\end{align} 
for $v\in V_{\pi_\fp}$ and  $d_\chiv= \ev_{d^\times a_\fp}^*(p_{(A_\fp,\chiv)})\in \Hom_{A_\fp}(V_\piv,\bc_\chiv)$. Note that while the integral \eqref{i-psi-chi-local-int} depends on the choice of the measure $dn_\fp$, the map \eqref{i-psi-chi-local} does not. The integral \eqref{i-psi-chi-local-int} is not absolutely convergent and should be understood in a regularized sense. 

We want to point out that the  integral \eqref{i-psi-chi-local-int} is not covered by the  Jacquet-Langlands  \cite{JL} theory, unlike the integral \eqref{i-chi-psi-local-int} for the opposite map.  
\end{remark}

\subsubsection{Global map} 

\begin{proposition*}\label{glob-i-map-prop-conv}  Fix an invariant rational differential form $\dl_{A}$ on $A$.
\begin{enumerate}
 \item   The tensor product $\fii(\psi,\chi, \dl_{A})=\otimes_\fp \fii(\psi_\fp,\chi_\fp, d^\times a_\fp)$ is absolutely convergent for $Re(\chi)\ll 1$, and has the meromorphic continuation to the complex space of all characters. 

\item The resulting map 
\begin{align*}\label{i-psi-chi-global}
\fii(\psi,\chi)=\fii(\psi,\chi, \dl_{A}): P(V_{\pi},L_\chi(A(\ba))) \to P(V_{\pi},L_\psi(N(\ba)))\ 
\end{align*}
does not depend on the choice of the rational form $\dl_{A}$. 
\end{enumerate}
\end{proposition*}

\subsubsection{Action on automorphic periods}
As the final step we have to compute the effect of the map $\fii(\psi,\chi)$ on the Hecke automorphic period. To this end, we now formulate our main result in this section.
\begin{theorem*}\label{T-per-N-per-thm} Let $\chi$ be such that the partial $L$-function satisfies $L_S(\pi,\chi)\not=0$ for some finite set $S$ of primes. The global map $\fii(\psi,\chi)$ is coherent, i.e., it sends the automorphic Hecke period $d_{\chi}$ to the automorphic Whittaker period $\CW^\psi$. Namely, we have
 \[ \fii(\psi,\chi)d_{\chi}=\CW^\psi\ .\]
\end{theorem*}

\begin{remarks}\label{a-tamag-rem} 1. We would like to point out  a subtle (but crucial) 
difference between Theorem~\ref{N-per-T-per-thm} and Theorem~\ref{T-per-N-per-thm}. The collection of local measures $\{d^\times a_\fp(\dl_{A})\}$ appearing in Proposition~\ref{glob-i-map-prop-conv} defines the Tamagawa adelic structure on the torsor $L^T(A(\ba))$ described in Section \ref{tam-meas-sect}, but these local measures do not define a {\it genuine  measure} on $A(\ba)$.  This differs from the situation described in Proposition \ref{glob-i-map-prop} where the Tamagawa adelic structure on $L^T(N(\ba))$ defines   a genuine measure on $N(\ba)$.
As a result, the direct map $\fii(\chi, \psi)$ from the Whittaker period space to the Hecke period space has the  integral representation (i.e., the  Hecke-Jacquet-Langlands integral \eqref{int} given by an integral over a closed (non-compact) cycle $A(k)\sm A(\ba)$ in the automorphic space $X$), but we do not know of such an integral representation for the map $\fii(\psi,\chi)$ in the opposite direction.  

2. Similarly to Remark \ref{rem-bundls-hol-sect} we can interpret Theorem \ref{T-per-N-per-thm} as a fiberwise meromorphic map from the  bundle  $H_\chi(\pi)$ to the constant bundle $W_{\psi,\chi}(\pi)$. This time however the corresponding map has poles. 
\end{remarks}
\subsubsection{The relation} We have the following relation between two maps we defined for Whittaker and Hecke periods.  
\begin{theorem*}\label{rel-thm} The following relation holds
\[ \fii(\psi,\chi)\circ \fii(\chi,\psi)=id\ , \]
as an endomorphism of $P(V_{\pi},L_{\psi}(N))$.
\end{theorem*}

\subsection{Proofs.} The logic of the proof we present is as follows.
We first prove results from the Section \ref{Wh-Hecke} by repeating arguments of Jacquet-Langlands in a slightly different language. We then prove Theorem \ref{rel-thm} by a local computation (see Lemma \ref{i*i=id-lemma}). This then implies all the other  results in Section~\ref{conv-Hecke-sect}.

\subsubsection{Regularization.}\label{reg-whitak-f} Many integrals appearing in our considerations require a regularization. There are various techniques to achieve regularization, e.g., via analytic continuation (see \cite{G1}).  In what follows we often use the following more elementary standard construction.

Let $F$ be a local non-archimedian field and $A$ the additive group of $F$. Fix a non-trivial  additive character $\psi: A \to \bc^\times$.

Let $V=V(A)$ denote the space of locally constant functions on $A$. It contains the subspace $\CS(A)$ of functions with compact support.

Let $V_0 \sub V$ be the space of functions $f \in V$ such that the image of $f$ in the space $V / \CS(A)$ is $F^\times$-finite.

Consider the Whittaker functional $W$ on the space $\CS(A)$ given by the integral $W(u)= \int_A \psi(a) \cdot u(a) da$.

\begin{claim}  There exists the unique extension of the Whittaker functional $W$ to the space $V_0$ that satisfies the identity
$W(t_a(u)) = \psi(a) W(u)$, where $t_a (u)(x) = u(x-a)$. 
\end{claim}

Uniqueness is easy since $F$ acts trivially on $V_0/\CS(A)$.  To prove existence one can extend the functional $W$ by setting  $W(u) = lim_{\{B\}}\int_B \psi(a) \cdot u(a) da$ where the limit is taken over a family $\{B\}$ of the standard compact balls exhausting $A$. It is easy to see that the limit exists for functions in $V_0$. 

\subsubsection{Proof of Proposition \ref{loc-i-map-prop}.} Both claims are standard in the Hecke-Jacquet-Langlands theory once the translation \eqref{i-chi-psi-local-int} into local zeta integrals is made.

\begin{enumerate}
	\item For a smooth vector $v\in V_{\pi_\fp}$, the Whittaker function $l^{\psi_\fp}(\piv(\bar a)v)$ is rapidly decreasing as $||t||\to\8$ in the positive Weyl chamber, and has a polynomial behavior in the opposite direction. This implies the absolute convergence of the integral for the character $\chi_\fp|\cdot|^{-s}_\fp$ with $Re(s)\gg 1$. 
	
The meromorphic continuation is equivalent to the meromorphic continuation of the Jacquet-Langlands local zeta integrals (see \cite{JL}). 
	
	\item This is the standard computation in the Hecke-Jacquet-Langlands theory. 
\end{enumerate}

\subsubsection{Proof of Proposition \ref{glob-i-map-prop}.} Indeed this follows immediately from the analytic continuation of $L(s,\pi)$ and from Proposition \ref{loc-i-map-prop}. In fact, this is a part of the Jacquet-Langlands method where the adelic integral is reduced to the absolutely convergent integral \eqref{int} via unfolding.

\subsubsection{Proof of Theorem \ref{N-per-T-per-thm}.} As we indicated before, the following proof is the standard Hecke-Jacquet-Langlands proof. On the basis of \eqref{p-T-def}, we want to compute
\begin{align}\label{1}
\left[i(\chi|\cdot|^{-s},\psi, dn,d^\times a)(\CW^\psi)\right](v)=\prod_\fp \int_{A_\fp}\chi_\fp\inv(\bar a)|a|_\fp^{s}\CW^{\psi_\fp}(\pi_\fp(\bar a)v)\ d^\times a_\fp \ ,
\end{align} for $Re(s)\gg 1$.  The absolute convergence of local integrals and the absolute convergence of the Euler product imply that
\begin{eqnarray}\label{11}
\prod_\fp \int_{A_\fp}\chi_\fp\inv(\bar a)|a|_\fp^{s}\CW^{\psi_\fp}(\pi_\fp(\bar a)v)\ d^\times a_\fp=\\ \nonumber \int_{A(\ba)}\chi\inv(\bar a)|a|^{s}\CW^{\psi}(\pi(\bar a)\phi_v)\ d^\times a \ .
\end{eqnarray}


We invoke the standard unfolding technique.
The rational torus acts transitively on Whittaker functionals for different characters. For an automorphic period $\CW^{\psi}$ and a character $\psi_\al(x)=\psi(\al x)$ with $\al\in k^\times$, we have the corresponding automorphic period given by $\CW^{\psi_\al}=\pi^{*}\left(
\left(\begin{smallmatrix}
\al & \\
                         &   1     \\
\end{smallmatrix} \right) \right)\CW^\psi$. We have  the following Fourier expansion at the identity for an automorphic function $\phi_v$, $v\in V_\pi$, in a cuspidal representation $\pi$, 
\begin{align}\label{Fourier-exp}
\phi_v(e)=\sum_{\al\in k^\times}\CW^{\psi_\al}(v)=\sum_{\al\in k^\times}\CW^{\psi}\left(\pi(\bar\al)v\right)\ .
\end{align}
 The  Fourier expansion for a cusp form $\phi_v$ implies that 
\begin{eqnarray}\label{2}
&&\int\limits_{A(\ba)}\chi\inv( \bar a)| a|^{s}\CW^{\psi}(\pi(\bar a)v)\ d^\times a =\\
&&\sum\limits_{\al\in k^\times}\int\limits_{A(k)\stm A(\ba)}\chi\inv(\bar a)|a|^{s}\CW^{\psi}\left(\pi(\bar\al )\pi(
\bar a )v\right)\ d^\times a =\nonumber\\
&&\int\limits_{A(k))\stm A(\ba)}\chi\inv(\bar a)|a|^{s}\left[\sum\limits_{\al\in k^\times}\CW^{\psi}\left(\pi(\bar\al )\left(\pi(
\bar a )v\right)\right)\right]\ d^\times a =\nonumber\\
&&\int\limits_{A(k)\stm A(\ba)}\chi\inv(\bar a)|a|^{s}\phi_{\pi(
\bar a )v}(e)\ d^\times a =\nonumber \\&& \int\limits_{A(k)\stm A(\ba)}\chi\inv(\bar a)|a|^{s}\phi_v(\bar a)\ d^\times a\ . \label{int}
\end{eqnarray} This gives the integral \eqref{torus-per-def} for the Hecke period
$d_{\chi|\cdot|^{-s}}$. 
We have used Fubini's theorem for $Re(s)\gg 1$, to decompose the adelic integral into the integral over a quotient space of the sum over $A(k)$  since all integrals are absolutely convergent. The resulting integral defines an analytic function for all values of $s$ and for a cuspidal $\pi$.

\subsubsection{Proof of Theorem \ref{T-per-N-per-thm}}
Proofs of all  statements leading to  and of the theorem itself are purely local, granted we already know the {\it direct} relation between Whittaker and Hecke periods (i.e., Theorem \ref{N-per-T-per-thm}). Namely, all proofs  follow from the following local lemma 
\begin{lemma*}\label{i*i=id-lemma} Let $\pi_\fp$ be an irreducible unitary representation  of $\GL_2(k_\fp)$, $\chi_\fp$ a quasi-character of  $A_\fp$,  $\psi_\fp$ a non-trivial character on $N_\fp$, and $dn_\fp$ and $d^\times a_\fp$ measures on $N_\fp$ and $A_\fp$, respectively, corresponding to some invariant differential forms $\dl_N$ and $\dl_A$. We have the following identity: 
\begin{align}\label{lem-identity} \fii(\psi_\fp,\chi_\fp,  d^\times a_\fp)\circ \fii(\chi_\fp,\psi_\fp,dn_\fp)= c(dn_\fp,d^\times a_\fp)\cdot id\ ,
\end{align} as an endomorphism of $P(V_{\pi_\fp},L_{\psi_\fp}(N_\fp))$.
Here  $c(dn_\fp,d^\times a_\fp)\in \bc$ is the proportionality  constant between the measure $d^\times a_\fp dn_\fp$ on $A_\fp\times N_\fp\simeq k^\times_\fp\times k_\fp$ and the measure $d^\times y dx$, i.e., $d^\times a_\fp dn_\fp= c(d^\times a_\fp, dn_\fp) d^\times ydx$. 
\end{lemma*}

\begin{remark}\label{adj-rem}
In this paper, we consider the construction of local maps $I_\fp$  by regularizing integrals over  appropriate subgroups (following the original construction of J.-L. Waldspurger \cite{Wa}). In certain cases one can construct such maps in both directions, i.e., maps $I_\fp(H_{1,\fp},H_{2,\fp}): \CP_\fp(\pi_\fp,\s_\fp)\to \CP_\fp(\pi_\fp,\tau_\fp)$ and $I_\fp(H_{2,\fp},H_{1,\fp}): \CP_\fp(\pi_\fp,\tau_\fp)\to \CP_\fp(\pi_\fp,\s_\fp)$. As was pointed out to us by Y. Sakellaridis, when $\tau_\fp$ and $\s_\fp$ are characters (as in examples in this paper), these maps are formally adjoint   in the following sense. We follow notations  from Section \ref{idea}. To an element $\xi\in \CP_\fp(\pi_\fp,\s_\fp)$, we can associate a map
$r_\xi:V_\piv\to \CF(H_{1,\fp}\setminus G)$ from the space of smooth vectors in the representation $\piv$ to the space of appropriate functions on $H_{1,\fp}\setminus G$ given by $v\mapsto \xi(\piv(g)v)$. Similarly for an element $\eta\in \CP_\fp(\pi_\fp,\tau_\fp)$, we have the map $q_\eta: V_\piv\to \CF(H_{2,\fp}\setminus G)$. The integration procedure $I_\fp(H_{1,\fp},H_{2,\fp})$ could be described then as the map $\CI_\fp(H_{1,\fp},H_{2,\fp}): \CF(H_{1,\fp}\setminus G) \to \CF(H_{2,\fp}\setminus G)$ given by the integral $[\CI_\fp(H_{1,\fp},H_{2,\fp})(\phi)](g)=\int_{H_{2,\fp}}\phi(hg)dh$ (we are leaving aside convergence issues). We have then $\langle \CI_\fp(H_{1,\fp},H_{2,\fp})(\phi),\psi\rangle_{H_{2,\fp}\setminus G_\fp}=\int\limits_{H_{2,\fp}\setminus G_\fp}\bigg(\int\limits_{H_{2,\fp}}\phi(hg)dh\bigg)\psi(g)dg=\int\limits_{G_\fp}\phi(g)\psi(g)dg= \int\limits_{H_{1,\fp}\setminus G_\fp}\phi(g)\bigg(\int\limits_{H_{1,\fp}}\psi(hg)dh\bigg)dg=\langle \phi,\CI_\fp(H_{2,\fp},H_{1,\fp})(\psi),\rangle_{H_{1,\fp}\setminus G_\fp}$.

Under certain conditions (which are satisfied for Whittaker/Hecke cases, we consider in Sections \ref{Hecke-sect} and \ref{conv-Hecke-sect}) which in \cite{SV} are called ``local unfolding'', it is shown by Y. Sakellaridis and A. Venkatesh \cite{SV} that above adjoint maps are also inverse of each other. In particular, this should imply our Lemma \ref{i*i=id-lemma} at least for tempered representations.  
\end{remark}
\subsubsection{Proof of Lemma \ref{i*i=id-lemma}}\label{PrfLem344} The proof follows from a direct computation in the Kirillov model (see Appendix \ref{Kir-mod-apx}).
In fact, this computation is essentially identical for all representations $\piv$ since it only involves the action of the Borel subgroup of $\GLt$. 

\subsubsection*{ Non-archimedian fields} Let $p_{N_\fp,\psiv}\in P(V_\piv, L_\psiv(N_\fp))$ and let $dn_\fp=dx$ be  the standard measure on $N_\fp\simeq  k_\fp$. The measure $dn_\fp$  induces the isomorphism $\ev_{dn_\fp}^*:P(V_\piv, L_\psiv(N_\fp))\to \Hom_{N_\fp}(V_\piv, \bc_\psiv)$. The Whittaker functional $\CW^\psiv=\ev_{dn_\fp}^*(p_{N_\fp,\psiv})$ gives rise to the Kirillov model realization $\Ck^{\CW^\psiv}:V_\piv\to \CK^\psiv(\piv)$ of $\piv$. Let $\dl_\fp=d^\times a_\fp=da_\fp/|a|_\fp$ be the standard local  measure  on $A_\fp\simeq k^\times_\fp$ and let $\ev_{\dl_\fp}^*:P(V_\piv, L_\chiv(A_\fp))\to \Hom_{A_\fp}(V_\piv, \bc_\psiv)$ be the corresponding isomorphism. 

We first compute the image $d^\#_{\chiv}$ of the Whittaker $\CW^\psiv$ functional under the integration with respect to $A_\fp$. We have $d^\#_{\chiv}:=(\ev_{\dl_\fp}^*)\inv(\fii(\chi_\fp,\psi_\fp, dn_\fp)\CW^\psiv)=\int_{A_\fp}\chi_\fp\inv(\bar a)\piv^*(\bar a)\CW^\psiv \dl_\fp$. In the Kirillov model of $\piv$, we have $\CW^\psiv(f)=f(1)$ for $f\in\CK^\psiv(\piv)$, and we have
\[ d^\#_{\chiv}(f)=\int_{A_\fp}\chi_\fp\inv(\bar a)\left[\piv(\bar a)f(x)\right]_{x=1} \dl_\fp=\int_{k_\fp^\times}\chi_\fp\inv(a)f( a) d^\times a_\fp \ .\]
Hence in the Kirillov model, the functional $(\ev_{\dl_\fp}^*)\inv(\fii(\chi_\fp,\psi_\fp, dn_\fp)\CW^\psiv)$ is given by the kernel $d^\#_{\chiv}(x)=\chi_\fp\inv(x)$ on $k^\times_\fp$. 

We now compute the image under the second integration with respect to $N_\fp$: $(\ev_{dn_\fp}^*)\inv( \fii(\psi_\fp,\chi_\fp, \dl_\fp)d^\#_{\chiv})$ given by 
\begin{align}\label{int1}\nonumber (\ev_{dn_\fp}^*)\inv( \fii(\psi_\fp,\chi_\fp, \dl_\fp)d^\#_{\chiv})(f)=&\int_{k_\fp}\psi_\fp\inv(x)d^\#_\chiv(\piv(n(x)f)dx=\\ &\int_{k_\fp}\psi_\fp\inv(x)\left[\int_{k_\fp^\times}\chi_\fp\inv(y)\psiv(xy)f(y)d^\times y\right]dx\ . 
\end{align} 
The inner integral is absolutely convergent for $Re(\chi_\fp\inv)\gg 1$, since functions in the Kirillov model are compactly supported on $k_\fp$ and have a polynomial behavior at $0$. 

We first compute the integral \eqref{int1} for functions which are compactly supported  on $k_\fp^\times$. Consider $f\in  \CS(k_\fp^\times)$ and assume that $f(u+y)=f(y)$ for all $y$ and  $|u|_\fp\leq q_\fp^{-N}$ for some $N\geq 0$. The inner integral in this case is zero for $|x|_\fp\geq q_\fp^{N+1}$. Hence we can take the outer integral over a big enough compact set $B_{N'}=\{|x|_\fp\leq q_\fp^{N'}\}$, $N'\geq N$. Both integrals are absolutely convergent over compact sets, and we can interchange their order.   We now have 
\begin{align}\label{int2}\nonumber (\ev_{dn_\fp}^*)\inv( \fii(\psi_\fp,\chi_\fp, \dl_\fp)d^\#_{\chiv})(f)=\int_{k_\fp^\times}\chi_\fp\inv(y)\left[\int_{B_{N'}}\psi_\fp((y-1)x)dx\right]f(y)d^\times y\ . 
\end{align} The inner integral is zero, unless $|y-1|_\fp\leq q_\fp^{-N'}$, and we have   $f(y)=f(1)$ under such a restriction. By possibly increasing $N'$, we also can assume that $\chi_\fp(y)\equiv 1$ for $|y-1|_\fp\leq q_\fp^{-N'}$. The integration with respect to measures $dx$ and $d^\times y$ then give $1$. Hence on the space  $\CS(k_\fp^\times)$, we have shown that 
\[(\ev_{dn_\fp}^*)\inv( \fii(\psi_\fp,\chi_\fp, \dl_\fp)d^\#_{\chiv})(f)=f(1)=\CW^\psiv(f)\ .\]
This finishes the proof of the Lemma for compactly supported functions (and finishes the proof for supercuspidal representations).
 
For induced representations, we also have to evaluate the integral \eqref{int1} on the space 
$V(\chi_1,\chi_2)$ describing the Kirillov model (see Section \ref{kir-mod-sect}). We claim that on the space $V(\chi_1,\chi_2)$ the integral \eqref{int1} is also given by the evaluation at the identity. Modulo compactly supported functions, functions spanning $V(\chi_1,\chi_2)$  are supported in $\CO_\fp$ and are essentially multiplicative characters near $0$. 
Hence we consider the integral 

\begin{align}\label{int3}\int_{k_\fp}\psi_\fp\inv(x)\left[\int_{\CO_\fp}\chi_\fp\inv(y)\psiv(xy)\chi'_\fp(y)d^\times y\right]dx\ ,
\end{align} for a fixed character $\chi_\fp'$. For $Re(\chi_\fp\inv)\gg 1$, the inner integral is absolutely convergent and decays polynomially in $|x|_\fp$. In fact, for $Re(\chi_\fp\inv\chi_\fp')\geq 2$, the inner integral  is bounded by $|x|_\fp^{-2}$ as $|x|_\fp\to\8$, and hence the outer integral is absolutely convergent. Hence for $Re(\chi_\fp\inv)\gg 1$, the functional defined by the integral \eqref{int1} extends to the space $V(\chi_1,\chi_2)$, and defines a functional on the whole space $\CW^\psiv(\piv)$. 
On the other hand, the double integral \eqref{int1} is clearly $(N_\fp, \psiv)$-equivariant. The space of such functionals is one-dimensional with a basis consisting of the Whittaker functional $\CW^\psiv$ we started with. Hence, for  for $Re(\chi_\fp\inv)\gg 1$, the integral \eqref{int1} coincides with the Whittaker functional $\CW^\psiv$ for induced representations as well. 

We  proved the statement \eqref{lem-identity} of the Lemma for $Re(\chi_\fp\inv)\gg 1$. Since the family of maps $ \fii(\chi_\fp,\psi_\fp, dn_\fp)$ is a meromorphic family of maps by \cite{JL}, the identity \eqref{lem-identity} holds for all $\chi_p$, and in fact provides the meromorphic continuation of the family of maps $ \fii(\psi_\fp, \chi_\fp, d^\times a_\fp)$.

Change of measures on $N_\fp$ and $A_\fp$ gives rise to the scaling factor $c(d^\times a_\fp,dn_\fp)$. 

\subsubsection*{Archimedian fields} We assume for simplicity that $k_\fp=\br$. The integral analogous to \eqref{int1} now takes the form
  \begin{align}\label{int1-arch}\int_{\br}e^{-ix}\left[\int_{\br_+}|y|^{-s}e^{ixy}f(y)d^\times y\right]dx\ ,
  \end{align} where $f$ is a smooth vector in the Kirillov space of a unitary irreducible representation of $\PGLR$. In particular, $f$ is smooth, rapidly decreasing as $y\to\8$, has a polynomial asymptotic at $0$ and $f\in L^2(\br_+, d^\times y)$. Hence, for $Re(s)\ll 0$, the inner integral is absolutely convergent and defines a bounded polynomially  decreasing function as $|x|\to\8$. This implies that the outer integral is also absolutely convergent for $Res(s)\ll 0$. Hence  \eqref{int1-arch} is given by the repeated Fourier transform and we obtain  $\CF\CF(f|\cdot|^{-s-1})(1)=f(1)$ for $Re(s)\ll 0$, i.e., the integral \eqref{int1-arch} gives the Whittaker functional in the Kirillov model. By analytic continuation we obtain this relation for all $s$. 

\section{Opposite Whittaker periods}\label{opp-whittak-sect} We denote by $N^+=\left\{
\left(\begin{smallmatrix}
1 & x\\
&    1      \\
\end{smallmatrix} \right)\right\}$ and  $N^-=\left\{
\left(\begin{smallmatrix}
1 & \\
x    &    1      \\
\end{smallmatrix} \right)\right\}$ two standard unipotent subgroups of $\GLt$. Our aim is to compare two period data coming from    Whittaker functionals associated with $(N^+,\psi)$ and $(N^-,\psi')$, where $\psi,\psi':\bq\sm \ba\to \bc$ are two non-trivial additive characters.

\subsection{Local map}\label{op-local-map} Let $\piv$ be an irreducible  representation of $G$, $\psi_\fp, \psi'_\fp:k_\fp\to\bc$ be two non-trivial additive characters. We denote by $\psi_\fp^+ :N^+_\fp\to\bc$ the character given by $\psi_\fp^+(n^+(x))=\psi_\fp(x)$ and by  $\psi_\fp^- :N^-_\fp\to\bc$ the character given by $\psi_\fp^-(n^-(x))=\psi'_\fp(x)$.  Consider the local Whittaker period space  $P^-_\fp({\pi_\fp},L_{\psi_\fp^-})=\Hom_{N^-_\fp}(V_{\pi_\fp},L_{\psi_\fp^-}(N_\fp^-))$. Choose an invariant differential form $\dl_-$ on $N^-$ and let  $dn^-_\fp=dn^-_\fp(\dl^-)$ be the corresponding invariant measure on $N^-$. We denote by  $\ev^*_{dn^-_\fp}:P^-(V_{\pi_\fp},L_{\psi_\fp^-}(N_\fp^-))\to \Hom_{N^-_\fp}(V_{\pi_\fp},\bc_{\psi_\fp^-})$ the induced isomorphism. We now  construct a map 
\begin{align}\label{ii-map} 
\fii(\psiv^+,\psiv^-,dn^-_\fp):  P^-(V_{\pi_\fp},L^-_{\psi_\fp^-})\to P^+(V_{\pi_\fp},L^+_{\psi_\fp^+})\ ,
\end{align} given by the integration.  Namely, for a vector $v\in V_\piv$ and a map $p_{\psiv}^-\in   P^-(V_{\pi_\fp},L^-_{\psi_\fp^-})$, we consider the following (matrix coefficient) function given by $\g_{p_{\psiv^-}^-,v}(n^+_\fp)= \ev^*_{dn^-_\fp}(p_{\psiv^-}^-(\piv(n^+_\fp)v))\in C^\8(N^+_\fp)$ and take its image in $L_{\psi_\fp^+}(N^+_\fp)$ under the integration/co-invariant map.

\begin{proposition*}\label{loc-W-W-prop}  \
\begin{enumerate}
	\item The map 
\begin{align}\label{w-map-prop} 
\fii(\psiv^+,\psiv^-, dn^-_\fp): P^-(V_{\pi_\fp},L^-_{\psi_\fp^-})\to P^+(V_{\pi_\fp},L^+_{\psi_\fp^+})\ ,
\end{align} is well-defined. 
\item 
For unramified $\piv$, $\psiv, \psiv'$ and $\dl_-$, we have 
\begin{align}\label{constant-fp-rel} \fii(\psiv^+,\psiv^-,dn^-_\fp)\xi^-_0 =\lm_\fp(\piv,\psiv,\psiv')\cdot \xi^+_0\ , \end{align} where $\xi^\pm_0$ are unramified and the constant 
$\lm_\fp(\piv,\psiv,\psiv')\in\bc$ is given by 

\begin{align}\label{constant-fp} \lm_\fp(\piv,\psiv,\psiv')=1-\tr(Ad(\s(\piv))q_\fp\inv\ , \end{align}
where $\s(\piv)\in\SL_2(\bc)$ denotes the Satake parameter of $\piv$ and $Ad$ is the adjoint representation of the dual group of $G$. 
\end{enumerate}
\end{proposition*}



\subsection{Global map} Let now $\pi=\otimes_\fp \piv$ be an automorphic cuspidal representation of $G$, $\psi=\otimes_\fp \psiv$ and $\psi'=\otimes_\fp \psiv'$ are two global non-trivial characters. We also choose a  non-zero invariant differential form $\dl^-$ on $N^-$. In order to construct the global map out of  local maps $\fii(\psiv^+,\psiv^-, dn^-_\fp)$, we need to glue constants $\{\lm_\fp\}$. The product $\prod_\fp\lm_\fp$ is not absolutely convergent. However, due to the unramified computation \eqref{constant-fp} one has the natural regularization procedure. This is based on the use of the adjoint $L$-function of $\pi$. For an unramified $\piv$, let $L(s,\piv,Ad)$ be the local adjoint $L$-function and  for  a finite set $S$ of primes, including primes where $\pi$ is ramified, let $L_S(s,\pi,Ad)$ be the partial adjoint $L$-function of $\pi$.

\begin{proposition*}\label{Euler-pr-thm}
 The Euler product
\begin{align}\label{Euler-pr-maps} \prod_{\fp\ {\rm unramified}}\left[{L(1,\piv,Ad)}\cdot\fii(\psiv^+,\psiv^-, dn^-_\fp)\right]
\end{align} is absolutely convergent. 
\end{proposition*}
On the basis of this proposition, we consider for a large enough set $S$, the following  absolutely convergent Euler product
\begin{align}\label{Euler-pr-map} \fii(\psi^+,\psi^-, \dl^-)=&&\\ L_S(1,\pi,Ad)\inv\prod_{\fp\in S}\fii(\psiv^+,\psiv^-, dn^-_\fp)&& \prod_{\fp\not\in S}\left[{L(1,\piv,Ad)}\cdot\fii(\psiv^+,\psiv^-, dn^-_\fp)\right]\ .
\nonumber\end{align} The resulting map clearly does not depend on a set $S$ if it is  large enough. We also use the well-known fact that $L_S(1,\pi,Ad)\not=0$.

\begin{theorem*}\label{N-N-indep-dl} 
The resulting map
\begin{align}\label{fii-map-glob} 
\fii(\psi^+,\psi^-)=\fii(\psi^+,\psi^-,\dl^-):  P^-(V_{\pi},L^-_{\psi^-})\to P^+(V_{\pi},L^+_{\psi^+})\ 
\end{align} is well-defined and
does not depend on  the choice of the form $\dl^-$. 
\end{theorem*}

\subsection{Action on automorphic periods. The invariant} We now consider the action of the map $\fii(\psi^+,\psi^-)$ on automorphic periods. According to the Theorem \ref{N-N-indep-dl}, there exists a constant $\lm(\pi,\psi,\psi')\in\bc$ such that \begin{align}\label{invariant} \fii(\psi^+,\psi^-)(\CW^{\psi^-})=\lm(\pi,\psi,\psi')\cdot\CW^{\psi^+} \ .
\end{align} We call it the {\it period invariant} associated to the map $\fii(\psi^+,\psi^-)$. The constant 
$\lm(\pi,\psi,\psi')$  gives rise to a  {\it global invariant} of $\pi$ (and of characters  $\psi$, $\psi'$) depending only on its automorphic realization (or even only on the isomorphism class of $\pi$ if it appears with  multiplicity one in the automorphic space, as is the case for $GL(2)$).

 The constant $\lm(\pi,\psi)$ measures to what extent the integration map $\fii(\psi^+,\psi^-)$ fails to be {\it coherent} (e.g., $\lm(\pi,\psi,\psi')=1$ if it is coherent, i.e., maps the automorphic period to the automorphic perid,  as was the case for maps considered in Section \ref{Hecke-sect}).  In Appendix \ref{Ramanj-apx}, we make a numerical evaluation of this invariant for the Ramanujan  holomorphic cusp form $\Dl$ of weight $12$ and level $1$. In particular, we will see that $\lm(\Dl,e^{2\pi ix})\not=1$, i.e., in that case the corresponding map is not {\it coherent}. This leads to a new numerical invariant of the Ramanujan function.

\subsubsection{Product formula}
We claim that the invariant $\lm(\pi,\psi,\psi')\in\bc$ could be computed via an absolutely convergent Euler product (i.e., it has local to global representation). Let $a\in k^\times$ be the scalar such that $\psi(x)=\psi'(ax)$. To write the product formula for  $\lm(\pi,\psi,\psi')$, we use the element 
$w_a=
\left(\begin{smallmatrix}
 & a\\
                  1      &        \\
\end{smallmatrix} \right)\in G$.  We know that $w_a$ maps the automorphic period $\CW^{\psi^-}$ to the automorphic period $\CW^{\psi^+}$ (i.e., $\pi^*(w_a)\CW^{\psi^-}=\CW^{\psi^+}$), and we also  know how $w_a$ acts on local period spaces.

Fix a local representation $\piv$ and  local character $\psiv, \psiv'$.
We have the  isomorphism of co-invariants $c(w):L_{\psi_\fp^+}(N^+_\fp)\simeq L_{\psi^-_\fp}(N_\fp^-)$ arising from the conjugation map $wN^+w\inv\to N^-$.  We have the natural map: 
\begin{align}\label{w-map} 
m(\psiv^+,\psiv^-): P^-(V_{\pi_\fp},L^-_{\psi^-_\fp})\to P^+(V_{\pi_\fp},L^+_{\psi^+_\fp})\ ,
\end{align} given by  action of the element  $w$, i.e.,   $\left[m(\psiv^+,\psiv^-)\xi\right](v)=c(w_a)(\xi (v))$  for any $v\in V_\piv$ and $\xi\in P^-(V_{\pi_\fp},L^-_{\psi^-_\fp}) $.
Hence there exists a constant $\lm_\fp(\piv,\psiv)\in\bc$ such that 
\begin{align}\label{lm-def-general}
\fii(\psiv^+,\psiv^-,dn^-_\fp)\xi= \lm_\fp(\piv,\psiv,\psiv')\cdot m(\psiv^+,\psiv^-)\xi\ ,
\end{align} for any $\xi\in P^-(V_{\pi_\fp},L^-_{\psi^-_\fp}) $. It is easy to see that for unramified $\fp$, these coefficients  coincide with  those defined in \eqref{constant-fp-rel}. 

\begin{theorem*}\label{Euler-pr-expan-thm}  
For a sufficiently large set $S$, the following relation holds  
\begin{align}\label{Euler-pr-S} \lm(\pi,\psi,\psi')=L_S(1,\pi,Ad)\inv\prod_{\fp\in S}\lm_\fp(\piv,\psiv,\psiv') \prod_{\fp\not\in S}\left[{\lm_\fp(\piv,\psiv,\psiv')}\cdot{L(1,\piv,Ad)}\right]\ ,
\end{align}
where $L_S(1,\pi,Ad)$ is the (analytically continued) partial adjoint $L$-function.

\end{theorem*}

\begin{remarks} 1. Instead of using the element $w_a$ to obtain the product formula for $\lm(\pi,\psi,\psi')$, it is possible to compare spaces $P^\pm(V_{\pi_\fp},L^\pm_{\psi^\pm_\fp})$ of Whittaker periods with  the space of  Hecke period (for some character) following the construction from Section \ref{Wh-Hecke}. This gives  another factoring of $\lm(\pi,\psi,\psi')$ into local factors. Such a comparison will differ locally from the one we considered above (i.e., coming from the action of $w_a$), and will  lead to different {\it local} factors. However,   {\it globally} this will give the same result due to the functional equation for the corresponding Hecke $L$-function of $\pi$. 
	
	2. The procedure of regularization described in Proposition \ref{Euler-pr-thm} is  similar to the one used to normalize the Tamagawa measure on the non-split torus (for the split torus the corresponding $L$-function has a pole) and is widely used in the theory of automorphic functions (e.g., see \cite{Wa}). 
	That the adjoint $L$-function shows up in our example is surprising since it is typically appear as a comparison between automorphic invariant Hermitian norm and the one coming from the  Whittaker model which is absent in our example. We discuss in Section \ref{aver-val-L} a possible geometric explanation for the appearance of the adjoint $L$-function. 
	
	 We have more examples of regularization of period maps similar to those appearing in Proposition~\ref{Euler-pr-thm}.  We hope to return to this subject elsewhere.

		3. Proposition \ref{Euler-pr-thm} could be formulated without mentioning $L$-functions explicitly, but using instead the language of maps between periods satisfying uniqueness property. The Rankin-Selberg method allows one to relate the adjoint $L$-function to the integration map from the diagonal period on $X_G\times X_G$ defining the invariant Hermitian form on $\pi\otimes \pi^\vee$ to the Hermitian form $\CW^\psi\otimes\overline\CW^\psi$ coming from the Whittaker functional. Hence the statement in Proposition \ref{Euler-pr-thm} could be interpreted as the statement about ratio for certain maps between appropriate  periods. The advantage of such a reformulation is that one does not need to know local components of $\pi$ in order to construct the regularization \eqref{Euler-pr-map} (i.e., one can think of the local factor $L(1,\pi_\fp,Ad)$ as a map between one-dimensional spaces of certain local period spaces). 
		
		4. As was pointed out by Y. Sakellaridis, the difference between examples in  Section \ref{Hecke-sect}  and the example form this section could be seen in the language of \cite{SV} as follows.  Following the general setup from Section \ref{idea}, we note that a choice of an invariant Hermitian form on $\piv$ (and on the relevant Gelfand data) gives rise to norm on the corresponding local period spaces $\CP_\fp(\pi_\fp,\s_\fp)$ and $\CP_\fp(\pi_\fp,\tau_\fp)$ (at least for tempered representations). Once the local map $I_\fp: \CP_\fp(\pi_\fp,\s_\fp)\to \CP_\fp(\pi_\fp,\tau_\fp)$ is constructed, one can ask if it is unitary with respect to these norms. It is easy to see that for the Hecke and converse to Hecke cases, the map is unitary and for the case of opposite Whittaker periods it is not unitary. 
		
		5. We can apply the construction of the global invariant \eqref{invariant} presented in this section  to an Eisenstein series representation $E(s)$ instead of a cuspidal representation $\pi$.   An argument similar to the proof presented above shows that the function $\lm(E(s),\psi,\psi')$ is well-defined for $s$ in the strip $|Re(s)|<\haf$ (in particular, $\lm(E(s),\psi,\psi')$ is well-defined for the unitary Eisenstein series, i.e., for $s\in i\br$). Hence we obtain the function $\lm(E(s),\psi,\psi')$ depending on the parameter $s\in\bc$.  We  note however that,  following the classical argument of T. Estermann \cite{E} (see also \cite{K1}, \cite{K2}),  one would expect that the Euler product \eqref{Euler-pr-S} $\prod_\fp\lm_\fp(\piv,\psiv,\psiv')$ with $\piv=\pi_s$ for all $\fp$, and unramified factors $\lm_\fp(\piv,\psiv,\psiv')$ as  in \eqref{constant-fp}, has the natural boundary at $Re(s)=\pm\haf$.
	\end{remarks}

\subsubsection{Functions on the dual group}\label{dual-gr-fn} Let us explain why it is natural to expect that the regularizing $L$-weights  should be given in terms of  Langlands $L$-functions. We stress again that we do not expect that for any period data as in Section \ref{local-step} some  regularizing $L$-weights exist.   

Instead of dealing with the local period data as in Section \ref{local-step}, we will work in  this section with the equivalent setup consisting of a group  $\mathfrak{G}=G_\fp\times {H}_{1,\fp}\times G_\fp\times {H}_{2,\fp}$, its subgroup $\mathfrak{H}={H}_{1,\fp}\times {H}_{2,\fp}\subset \mathfrak{G}$ under the natural embedding and the representation  $\Pi_\fp=\piv^*\otimes\chi_{1,\fp}\otimes\piv\otimes\chi_{2,\fp}^*$ of  $\mathfrak{G}$. We will assume for simplicity that $\fG$ is a split reductive group (we note that some of our examples do not fit this restriction, but one can easily modify this setup). We have then $\CP_\fp(\Pi_\fp,\bc)\simeq\CP_\fp(\pi_\fp,\chi_{1,\fp} )^*\otimes \CP_\fp(\pi_\fp,\chi_{2,\fp} )$. The local procedure of Section  \ref{local-step} should provide us with the canonical element $\CI_\fp\in \CP_\fp(\Pi_\fp,\bc)$ corresponding to the canonical map $I_\fp: \CP_\fp(\pi_\fp,\chi_{1,\fp} )\to  \CP_\fp(\pi_\fp,\chi_{2,\fp})$ and resulting local  unramified factors $D_\fp(\Pi_\fp)$ as in \ref{reg-by-L}. Our scheme calls for a regularization of  the Euler  product $\prod_{\fp\not\in S} D_\fp(\Pi_\fp)$ by the Langlands $L$-functions. The main issue here is how to identify the correct local $L$-function. Note that we do not assume that there is a natural parameter $s\in\bc$ present in the problem and this complicates the matter. Sometimes there is such a parameter (e.g., the parameter $s$ of  the character $|\cdot|^s$ for a split tori) and in that case one can avoid the discussion below of the extended dual group. 

We view  local Langlands $L$-functions as functions on the extended dual group $ \tilde\fG=\check\fG\times \br^\times_+$, where $\check\fG$ is the dual group of $\fG$.  For  a finite dimensional complex representation $r: \check\fG\to End_\bc(V_r)$ and $s\in\bc$, we consider a function on 
$\tilde\fG$ given by $L(s,\check g\times t,r)=\left[\det(1-r(\check g))t^{-s})\right]\inv$.
For an unramified representation $\Pi_\fp$ of $\fG_\fp$, the local $L$-factor is then given by $L_\fp(s,\Pi_\fp,r)=L(s,\s(\Pi_\fp)\times q_\fp^{-s},r)$, where $\s(\Pi_\fp)\in \check\fG$ is the Satake parameter of $\Pi_\fp$ and $q_\fp$ is the size of the residue field at $\fp$. 

We denote by $\mathcal{C}(\tilde \fG)$ the space of complex valued regular central functions on $\tilde \fG$.
Functions $L(s,\check g\times t,r)$ form a basis of  $\mathcal{C}(\tilde \fG)$.

Let us assume that local factors $D_\fp(\Pi_\fp)$ are defined for a rich family of unramified representations $\Pi_\fp$ (e.g., for all unramified tempered representations). Under such an assumption, for each $\fp$, we can view $D_\fp(\Pi_\fp)$ as a function defined on an open set $\Sigma_\fp\subset \check\fG$  of the dual  group 
and resulting functions are central. Hence we obtain a collection of central functions $\{D_\fp\}_{\fp\not\in S}$ each defined on the corresponding sheet of the set  $ \Sigma=\bigcup_{\fp\not\in S}\left[\Sigma_\fp\times q_\fp\right] \subset \tilde \fG$. Any regular central function on $\tilde\fG$ is uniquely determined by its values on the set $\Sigma$.

In view of J. Bernstein rationality theorem (see  \cite{Be}, \cite{Ban}), it is natural to expect that function $D_\fp(\Pi_\fp)$ is  a rational function in $q_\fp^{-s}$, and such a function is determined by its values 
for an infinite set of unramified $\fp$ and a fixed $s\not=0$. 

For a given local factors $D_\fp(\Pi_\fp)$ coming from an unramified computation in the local step \ref{local-step}, it is easy to check if the factorization into local Langlands $L$-functions  exists since a rational central function on $\tilde\fG$ is uniquely determined  by its values on {\it infinite} number of sets  $\Sigma_\fp\times q_\fp\subset \tilde\fG$ with $\fp$ unramified. The use of Langlands local factors  is natural since functions $L(s,\check g\times t,r)$ form a basis for central functions.

The above discussion supports our expectation that  Langlands $L$-functions could provide regularization in our global step, at least in some cases. This does not however explicitly describe the relevant $L$-factor. In Section \ref{aver-val-L} we will discuss an example where we can identify the regularizing factor by a geometric construction.

\subsubsection{Geometry of periods and the regularizing factor}\label{aver-val-L}
We present a geometric observation of how one can see what $L$-factor appears naturally in the regularization for two Whittaker periods   discussed above. We hope to discuss more examples elsewhere. 

In the local step \ref{local-step} we construct a canonical map $I_\fp: \CP_\fp(\pi_\fp,\chi_{1,\fp})\to \CP_\fp(\pi_\fp,\chi_{2,\fp})$. In particular, we obtain the collection of unramified maps $d_\fp\in \CP_\fp(\pi_\fp,\chi_{2,\fp})$ defined for  unramified primes $\fp$. The local regularizing factor then comes from this unramified computation. We note that in fact it is natural to ask how this factor varies as we change the period data. Namely, pairs of period data $(G,\pi,{H_1},\chi_1)$
and $(G,\pi,{H_2},\chi_2)$ naturally form  a ``moduli" space (or a stack). 

For the Hecke case from Section \ref{Hecke-sect},  we have a pair of a tori and a unipotent subgroup in the same Borel subgroup (and their corresponding characters). It is easy to see that this moduli space consists of a point and the regularizing factor will not change as we vary  the tori in an {\it unramified} fashion (and change accordingly the character), i.e., conjugate the tori by an element in $N_\fp\cap K_\fp$.    On the other hand, in Section  \ref{opp-whittak-sect} we consider two  Whittaker functionals for unipotent subgroups which are not in the same Borel (in the example we choose them to be conjugated by the Weyl element $w_a$). In this case the moduli space is non-trivial. If we fix one such a functional, the second one is given up to a conjugation by an element  $g\in N_\fp\sm G_\fp$. To consider an unramified data, we should conjugate by an element  $g\in (N_\fp\cap K_\fp) \sm K_\fp$. This correspond to the compact part of the moduli space (i.e., its points over $\bz_\fp$). It turns out that unlike in the Hecke case, here the result of an unramified computation depends on the conjugating element $g$ (i.e., it depends on the relative position of two unramified Whittaker functionals). We can average the resulting factor $d_\fp$ over the compact part of the moduli space, i.e., over $(N_\fp\cap K_\fp)\sm K_\fp$. As in Lemma \ref{i-to-j-funct1}, $d_\fp$ is given by the Bessel function at the point $g\in G_\fp$ depending on two Whittaker functionals.  The Bessel function at $g$ represents  essentially the (regularized) value of the  pairing $\langle \pi_\fp(g) \CW^{\psi},\CW^{\bar\psi}\rangle$.  Hence,  the averaged over  $g\in (N_\fp\cap K_\fp)\sm K_\fp$ value  is given by the absolute value squared of a Whittaker functional evaluated at the standard unramified vector. This gives the factor $L(1,\piv,Ad)$ according to the well known computation (see, e.g. \cite{Jacq01}) which is exactly the regularizing $L$-factor from Theorem \ref{Euler-pr-expan-thm} and might be the reason we see the adjoint $L$-function in \eqref{Euler-pr-S}.

\subsection{Proofs} Proofs below are based on repeated use of regularization of Whittaker functional as formulated in Section \ref{reg-whitak-f}.

\subsubsection{Local maps} Functions $\g_{p_{\bar\psiv}^-,v}$ defined in Section \ref{op-local-map} are not compactly supported on $N^+_\fp$ and corresponding integrals  should be understood in the regularized sense as in Section \ref{reg-whitak-f} . 

Denote by $\dl_+=w^*_a\dl_-$  the form on $N_+$ and let $dn^+_\fp$ be the corresponding invariant measure on $N^+_\fp$. We denote by $\ev^*_{dn^+_\fp}:P^+(V_{\pi_\fp},L^+_{\psi_\fp})\to \Hom_{N^+_\fp}(V_{\pi_\fp},\bc_{\psi_\fp})$ the induced isomorphism.  Hence we obtain the map 
\begin{align}\label{i-map} 
i(\psiv^+,\psiv^-,dn^-_\fp, dn_\fp^+): \Hom_{N^-_\fp}(V_{\pi_\fp},\bc_{\psi'_\fp})\to \Hom_{N^+_\fp}(V_{\pi_\fp},\bc_{\psi_\fp})\ ,
\end{align} given by $i(\psiv^+,\bar\psiv^-,dn^-_\fp, dn_\fp^+)=
\ev^*_{dn^+_\fp}(\fii(\psiv^+,\psiv^-,dn^-_\fp))$ which has the following integral representation 
\begin{align}\label{int-W-W} \left[i(\psiv^+,\psiv^-,dn^-_\fp, dn_\fp^+)p^-_{\psiv}\right](v)=&\int_{N^+_\fp}\bar{\psi}_\fp(n^+_\fp)\g_{p_{\bar\psiv}^-,v}(n^+_\fp)dn^+_\fp\ .
\end{align} The integral \eqref{int-W-W} does not converge absolutely, and  should be  regularized. For a non-archimedian field $k_\fp$, we will understand under the integral \eqref{int-W-W} the limit 
\begin{align}\label{int-W-W-B-N}\lim\limits_{l\to\8}\int_{B_l}\bar{\psi}_\fp(n^+_\fp)\g_{p_{\bar\psiv}^-,v}(n^+_\fp)dn^+_\fp\ ,
\end{align} where
$B_l=\{n(x),\ |x|_\fp\leq q_\fp^{l}\}\subset N^+_\fp$. We will show that for any given smooth vector $v\in V_{\pi_\fp}$, the integral \eqref{int-W-W-B-N} stabilizes as $l\to\8$. 

For $k_\fp=\br$, we will use the asymptotic expansion of the Bessel function of $\pi_\8$ in order to regularize integral \eqref{int-W-W} by means of analytic continuation (this procedure could be interpreted as the analytic continuation of the integral \eqref{int-W-W} in the space of parameters of representations of  $\GLR$).

\subsubsection{Proof of Proposition \ref{loc-W-W-prop}}\label{proof-N-N-1} The proof of the proposition is based on the same idea as the proof of Lemma \ref{lem-identity}, i.e., we compute the local map in terms of the Bessel function.

\begin{lemma*}  The following relation holds
	\begin{align}\label{i-to-j-funct1}
\fii(\psiv^+,\psiv^-,dn^-_\fp)\xi= j_{\piv,\psiv}(a)\cdot m(\psiv^+,\psiv^-)\xi\ ,
\end{align} for any $\xi\in P^-(V_{\pi_\fp},L^-_{\psi_\fp^-})$. Here $j_{\piv,\psiv}$ is the Bessel function of the representation $\piv$ (see Appendix \ref{Kir-mod-apx}) and $ m(\psiv^+,\psiv^-)$ is given by the action of $w$ as in \eqref{w-map}. 
\end{lemma*}

The lemma clearly implies that the local map is well-defined. From the lemma we see that $\lm_\fp(\piv,\psiv,\psiv')= j_{\piv,\psiv}(a)$. 

\subsection*{ Proof of the lemma.} 
\subsubsection{ Non-archimedian field} We fix an additive character $\psiv$ and choose a non-zero $(N_\fp^+,\psiv)$-Whittaker functional $\CW_+^\psiv$ on $\piv$. Let $\CK^\psiv (\piv)$ be the corresponding $(N_\fp^+,\psiv)$-Kirillov model of $\piv$. In the Kirillov model, the original Whittaker functional is then given by the delta function $\dl_1$ at $1\in k_\fp^\times$.  Let $j_{\piv,\psiv}$ be the $\psiv$-Bessel function of $\piv$. We have 
\begin{align}\nonumber \left[m(\psiv^+,\psiv^-)\dl_1\right](v):=\pi^*_\fp(w_a)\dl_1(v)=\dl_1(\piv(w_a)v)\\ =\langle \dl_{\al=a}, \int_{k^\times_\fp}j_{\piv,\psiv}(\al t)v(t)d^\times t\rangle\ .
\end{align} Hence the functional $\dl_-=m(\psiv^+,\psiv^-)\dl_1$ is given by the kernel $j_{\piv,\psiv}(at)$ in the $(N_\fp^+,\psiv)$-Kirillov model.
We now compute
\begin{align}\label{int-W-W-2}\nonumber \left[\fii(\psiv^+,\psiv^-)\dl_-\right](v)=&\int_{N^+_\fp}\bar{\psi}_\fp(x)\dl_-(\pi_\fp( n_+(x))v)dx\\ &=\int_{k_\fp}\bar\psiv(x)\left[\int_{k^\times_\fp}\psiv(xt)j_\piv(at)v(t)d^\times t\right]dx\ .
\end{align} For $v\in\CS(k^\times_\fp)$, this immediately implies that $\fii(\psiv^+,\psiv^-)\dl_-=j_\piv(a)\dl_1$ as in the proof of Lemma \ref{i*i=id-lemma}. For a non-archimedian field,  this finishes the proof for the space of compactly supported functions.  

For induced representations over a non-archimedian field, we note that the inner integral in \eqref{int-W-W-2} is absolutely convergent. This follows from the bound \eqref{j-bound} $|j_\piv(t)|\leq C_\piv|t|_\fp^{-1/4}$  (we assume that the central character is trivial), and the fact that functions in $V(\chi_1,\chi_2)\subset \CK^\psiv(\piv(\chi_1,\chi_2))$ satisfy the bound $|f(t)|\ll|t|_\fp^{1/2}\log|t|_\fp$ (both bounds hold for small enough $|t|_\fp$). Hence we can interchange the order of integration in \eqref{int-W-W-2} if we understand under the outer integral the limit $\lim_{N\to\8}\int_{|x|_\fp\leq p^N}\dots$\ . For every $N\geq 1$, we consider the absolutely convergent double integral 
\begin{align}\nonumber \int_{|x|_\fp\leq q_\fp^N}\int_{k^\times_\fp}\bar\psiv(x)\psiv(xt)j_\piv(at)v(t)d^\times t dx\ .
\end{align} Integrating now over $x$ first, we see that for any given smooth function $v$, the integral stabilizes as $N\to\8$. 
Hence the functional $\fii(\psiv^+,\psiv^-)\dl_-$ extends to the space $\CK^\psiv(\piv)$ for induced  representations as well. The uniqueness of the Whittaker functional implies again that the resulting functional is $\dl_1$, and hence we 
proved that for any unitary infinite-dimensional representation $\piv$ of $G$ over a non-archimedian field, the following relation holds:
\begin{align}\label{i-to-j-funct2}
\fii(\psiv^+,\psiv^-)\circ m(\bar\psiv^-,\psiv^+)\dl_1=\fii(\psiv^+,\psiv^-)\dl_-=j_\piv(a)\dl_1\ .
\end{align}

\subsubsection{Archimedian field} We now prove the same statement over reals. As before we have $\dl_-=\pi(w_a)\dl_1=j_{\pi_\8,\psi_\8}(a\cdot)$, and we consider the integral 
\begin{align}\label{int-W-W-arch}\nonumber \left[\fii(\psi_\8^+,\psi_\8^-)\dl_-\right](v)=&\int_{N^+_\8}\bar{\psi}_\8(x)\dl_-(\pi_\8(n_+(x))v)dx=\\ &\int_{\br}\bar\psi_\8(x)\left[\int_{\br^\times}\psi_\8(xt)j_{\pi_\8}(at)v(t)d^\times t\right]dx\ .
\end{align}
For Schwartz functions $v\in\CS(\br^\times)$, the integral is absolutely convergent 
and rapidly decaying in $|t|\to\8$. Hence we can split the outer integral into a compact part $|t|\leq N$ and the rest: $|t|>N$. The non-compact part tends to $0$ as $N\to\8$, and in the compact part we can change the order of integration. As a result, we arrive at $\left[\fii(\psi_\8^+,\psi_\8^-)\dl_-\right](v)=j_{\pi_\8,\psi_\8}(a) \dl_1(v)$ as in the non-archimedian case.  We need
to show that $i(\psi_\8^+,\psi_\8^-)\dl_-$ extends to a functional on $\CK^{\psi_\8}(\pi_\8)$. The inner integral in \eqref{int-W-W-arch} is absolutely convergent for all $v\in \CK^{\psi_\8}(\pi_\8)$ as follows from asymptotic of Whittaker functions of smooth vectors and from asymptotic of Bessel functions (e.g., asymptotic \eqref{j-asymp} for the $J$-Bessel function). Using these asymptotic, we see that the inner integral also has polynomial asymptotic expansion of the type $\sum_{i=0}^M a_i|t|^{\lm-i}+O(|t|^{Re(\lm)-M-1})$ as $|t|\to\8$ where $\lm\in\bc$ is the parameter of the representation $\pi_\8$. Such an integral could be regularized  by the analytic continuation method (see \cite{G1}). Hence we extended the functional  to the whole space  $\CK^{\psi_\8}(\pi_\8)$, and from the uniqueness of Whittaker functional it follows that  $\fii(\psi_\8^+,\bar\psi_\8^-)\dl_-=j_{\pi_\8,\psi_\8}(a) \dl_1$. \qed

\begin{remark}One can use the asymptotic  expansion for the Bessel function obtained in \cite{JY}, Proposition 2.3,  to give a proof for the Lemma for a non-archimedian field arguing as in the case of reals. 
\end{remark}

\subsubsection{Proof of \eqref{constant-fp-rel}}  This is a simple computation following \cite{S}, \cite{BM}. The Bessel function of an induced representation $\piv(\chi_1,\chi_2)$ is given by
\begin{align}\label{j-1-int} j_\piv(1)=\lim_{N\to\8}\int_{p^{-N}\leq|t|_\fp\leq p^N}\chi_1\inv\chi_2(t){\psi}_\fp(x-x\inv)|t|_\fp\inv dt\ .
\end{align} The integral stabilizes as $N\to\8$, and 
for unramified $\psiv$ and $\chi_i$, in fact, stabilizes at $N=1$.  We have $j_\piv(1)=\int_{p^{-1}\leq|t|_\fp\leq p}\chi_1\inv\chi_2(t){\psi}_\fp(x-x\inv)|t|_\fp\inv dt=\int_{|t|_\fp= 1}1 dt+\chi_1\inv\chi_2(\varpi)\int_{|t|_\fp= p\inv}{\psi}_\fp(-x\inv)|t|_\fp\inv dt+\chi_1\inv\chi_2(\varpi\inv)\int_{|t|_\fp= p}{\psi}_\fp(x)|t|_\fp\inv dt$. We obtain correspondingly 
\begin{align}\label{j-1-unram} j_{\piv,\psiv}(1)=1-q_\fp\inv-\chi_1\inv\chi_2(\varpi)q_\fp\inv-\chi_1\chi_2\inv(\varpi)q_\fp\inv\ .
\end{align}

\subsubsection{Proof of Proposition \ref{Euler-pr-thm}}
Let $\piv\simeq\piv(\chi_1,\chi_2)$ be an unramified representation. Denote by $\al^2_\fp=\chi_1\inv\chi_2(\varpi)$. We have \[ L(1,\piv,Ad)=1/(1-q_\fp\inv)(1-\al^2_\fp q_\fp\inv)(1-\al_\fp^{-2} q_\fp\inv)\ .\] Using  \eqref{constant-fp}  we write $\lm_\fp(\piv,\psiv)=1-q_\fp\inv-\al_\fp^2q_\fp\inv- \al_\fp^{-2}q_\fp\inv=L(1,\piv,Ad)\inv-r_\fp(\al_\fp, q_\fp)$, where $r_\fp(\al_\fp, q_\fp)=q_\fp^{-2}+\al_\fp^2q_\fp^{-2}+\al_\fp^{-2}q_\fp^{-2}-q_\fp^{-3}$. From this we deduce that 
\begin{align}\label{lm-to-L}&\lm_\fp(\piv)L(1,\piv,Ad)=(L(1,\piv, Ad)\inv-r_\fp(\al_\fp, q_\fp))L(1,\piv,Ad)=\\ & 1-r_\fp(\al_\fp, q_\fp)L(1,\piv,Ad)=1-Q_\fp(\al_\fp,q_\fp)\ . \nonumber
\end{align} Note that the term $Q_\fp(\al_\fp,q_\fp)=q_\fp^{-2}\frac{(1+\al_\fp^2+\al_\fp^{-2}-q_\fp^{-1})}{(1-q_\fp\inv)(1-\al^2_\fp q_\fp\inv)(1-\al_\fp^{-2} q_\fp\inv)}$ is expected to be bounded by  $q_\fp^{-2+\eps}$ according to the Ramanujan-Peterson  conjecture, and hence this leads to an absolutely convergent Euler product. Namely,  according to any non-trivial bound towards Ramanujan, there exists $\s>0$ such that $|\al_\fp|\leq q_\fp^{\haf-\s}$. Hence $|Q_\fp(\al_\fp, q_\fp)|\ll q_\fp^{-1-\s'}$ for some $\s'>0$. This implies that $|\lm_\fp(\piv)L(1,\piv,Ad)|\leq 1-q_\fp^{-1-\s'}$ and hence the Euler product $\prod_\fp \lm_\fp(\piv)L(1,\piv,Ad)$ is absolutely convergent.

\appendix
\section{Computation for the Ramanujan cusp form}\label{Ramanj-apx}

\subsection{Numerical evaluation}\label{Ram} Let $\Delta(z)=\sum_{n\geq 1}\tau(n)q^n$ be the classical  cusp form, studied by Ramanujan in \cite{Ra}, with $\tau(n)$ the Ramanujan tau function. The holomorphic cusp form $\Dl$ has weight $12$ and level $1$.  In the adelic language, $\Dl$ corresponds to a cuspidal automorphic representation $\pi_\Dl=\otimes_{p\leq \8}\pi_p$ of $GL(2)$ over $\bq$ with trivial central character. The corresponding local components are unramified for all finite primes $p$, and $\pi_\8$ is isomorphic to the discrete series representation of $\GLR$ with the lowest weight vector of weight $12$ (in \cite{BM} such a representation is denoted  by $\pi_6$ ; see Appendix \ref{real-Kir-apx} below). The Satake parameters $(\al_p, \al_p\inv)$ of a local representation $\pi_p$ at $p< \8$ are given by $\al_p+\al_p\inv=\tau(p)p^{-\frac{11}{2}}$ (since the Ramanujan conjecture is known for $\Dl$, we have $|\al_p+\al_p\inv|\leq 2$). Below we attempt to calculate the constant $\lm(\pi_\Dl,\psi)$ for the additive character $\psi(x)=e^{2\pi ix}$ of $\bq\stm \ba_\bq$.   We have 
\begin{align*}
\lm_p(\pi_p,\psi_p)=1-p\inv-\al_p^2 p\inv- \al_p^{-2} p\inv= 1-(1+\al_p^2+\al_p^{-2})p\inv\\ =1-((\al_p+\al_p\inv)^2-1)p\inv=1-(\tau^2(p)p^{-11}-1)p\inv=1-\tau(p^2)p^{-12}\ .
\end{align*} We also have 
\begin{align*}&L(1,\pi_p, Ad)=1/(1-p\inv)(1-\al_p^2 p\inv)(1-\al_p^{-2}p\inv)=\\ &(1-p\inv-\al_p^2 p\inv- \al_p^{-2} p\inv+ \al_p^2 p^{-2}+\al_p^{-2} p^{-2}+p^{-2}-p^{-3})\inv=\\ &(1-(\tau^2(p)p^{-11}-1)p\inv+(\tau^2(p)p^{-11}-1)p^{-2}-p^{-3})\inv \ .
\end{align*}  

We now consider a finite  product over first $N$ primes 
\begin{align}
\tilde\lm^f_N(\pi_\Dl,\psi)=\prod_{i\leq N}\lm(\pi_{p_i},\psi_{p_i})L(1,\pi_{p_i},Ad)
\end{align} for $N\geq 1$, and $\lm^f_N(\pi_\Dl,\psi)=\tilde \lm^f_N(\pi_\Dl,\psi)/L(1,\pi,Ad)$.  We have $L(1,\pi,Ad)= 0.63179294573\dots$ (see \cite{Za} and  Rubinstein's Sage $L$-functions Calculator). We computed a numerical approximation  $\lm^f_{100}(\pi_\Dl,\psi)=1.49154\dots$, and the archimedian counterpart $\lm^\8(\pi_{6}, \psi_\8)$ given by the value of the classical Bessel function  $j_{\pi_d,\psi_\8}(1)=2\pi J_{11}(4\pi)=1.8305\dots\ $ (see \eqref{j-disc}).

Hence we obtain the following numerical approximation for $\lm(\Dl, e^{2\pi ix})=4.32\dots\ $. In particular, one can see that it is different from $1$ (since it is easy to estimate the absolutely convergent remainder). 

\subsection{An infinite product}\label{strange} We do not understand the nature of the constant $\lm(\pi,\psi)$. According to the local unramified computation (see \eqref{constant-fp}),
the local constant $\lm_\fp (\piv,\psiv)$ and the Euler polynomial $L(1,\piv,Ad)$ have the same linear part. This allowed us to prove Proposition \ref{Euler-pr-thm} by approximating  $\lm_\fp(\piv,\psiv)$ with $L(1,\piv, Ad)\inv$. We can iterate this process. 

We consider the polynomial $l(a,x)=1-(a+1+a\inv)x$ in $x$ (so that $\lm_\fp(\pi_\fp(\al_\fp,\al_\fp\inv))=l(\al_\fp^2,q_\fp\inv$)). Let us introduce a family of polynomials 
\begin{align} p_l(a,x)=\prod_{i=-l}^l
(1-a^ix)\ ,\ l\geq 0
\end{align} (e.g., $p_0(a,x)=1-x$, $p_l(a,x)=(1-a^lx)(1-a^{-l}x)p_{l-1}(a,x)$). In particular, we have $L(s,\pi_\fp(\al_\fp,\al_\fp\inv), Sym^{2l})=p_l(\al_\fp^2,q_\fp^{-s})\inv$ for the symmetric power $L$-function. From an easy inductive argument, it follows that there are integer coefficients  $m_{kl}\in\bz$, $k,\ l\in \bz_+$ such that
\begin{align}\label{comb-iden-polyn}
l(a,x)=p_1(a,x)\prod_{k=2}^\8\left[\prod_{l=0}^{k-1}p_l(a,x^k)^{m_{kl}}\right] \ 
\end{align} as a formal power series identity in $\bz[a+a\inv][[x]]$. 

We now introduce  the constant $\lm^f(\pi,\psi)=\prod\limits_{\fp<\8}\lm_\fp(\piv,\psiv)$. 
Assuming coefficients $m_{kl}$ do not grow too fast (although they do grow exponentially), the formula \eqref{comb-iden-polyn}  suggests the following
representation (at least in the unramified case):
\begin{align}\label{comb-iden}
\lm^f(\pi,\psi)=L(1,\pi, Sym^{2})\prod_{k=2}^\8\left[\prod_{l=0}^{k-1}L(k,\pi, Sym^{2l})^{-m_{kl}}\right] \ .
\end{align} We note that all $L$-functions appearing in the infinite product are in the region of the absolute convergence (assuming the Ramanujan conjecture).

The sequence $m_{kl}$  could be interpreted as a virtual representation of $\SL(2,\bc)\times\bg_m$.
Even some basic properties of the sequence $m_{kl}$ are not clear to us. In particular, we do not know if  these coefficients are non-negative (i.e., is it true that $m_{kl}\geq 0$; this would mean that the corresponding virtual representation is a genuine representation). Also we would like to have an estimate for  the growth rate of $m_{kl}$ in order to justify convergence of the infinite product \eqref{comb-iden}. 

Here we list the first few coefficients $m_{kl}$ (kindly computed by S. D. Miller): \vskip 10pt

{\tiny
\begin{center}
  \begin{tabular}{ c | c c c c c c c c c c c }
    \hline
    $k\stm l$ & 0&1&2&3&4&5&6&7&8&9&10 \\ \hline
    1 &   & 1 \\ 
    2 &   & 1 \\
		3& & 1&1\\
		4& &2&1&1\\
		5& 1&3&3&2&1\\
		6&1&7&6&5&2&1\\
		7&5&13&15&12&7&3&1\\
		8 &9& 31& 33& 31& 18 &10& 3& 1\\
		9 &25& 67& 84& 74& 52& 29& 12& 4& 1\\
10 &55& 163& 198& 192& 137& 85& 39& 16& 4 &1\\
11 &144& 383& 500& 483& 375& 240& 127& 55& 19& 5& 1\\
  \end{tabular}
\end{center}
}

\section{Kirillov model }\label{Kir-mod-apx}

  Here we collect various facts about the Kirillov model for representations of $\GLt$ (for more detail, see \cite{JL}, \cite{Be}, \cite{Ba}, \cite{BS}, \cite{BM}, \cite{BM-real}, \cite{S}). 

\subsection{Non-archimedian Kirillov model}\label{kir-mod-sect} Let $\pi_\fp$ be an irreducible infinite dimensional unitary representation of $\GLt$ over a local field $k_\fp$. We fix a non-trivial character $\psi_\fp$ of $N_\fp$, and choose a non-zero Whittaker functional $\CW^{\psi_\fp}$ on $\pi_\fp$. Such a functional gives rise to the Kirillov model for $\pi_\fp$. Let $ S^+(k_\fp^\times)$ be the space of smooth (locally constant for $\fp<\8$) functions on $k_\fp^\times$ of rapid decay at infinity (relative to the completion $k_\fp^\times\subset k_\fp$ at $0$).  Consider the map $\Ck^{\CW^\psiv}: V_\piv \to S^+(k_\fp^\times)$ given by $\left(\Ck^{\CW^\psiv}(v)\right)(a)=\CW^\psiv(\piv(\bar a))$, $\bar a=
\left(\begin{smallmatrix}
a & \\
                         &    1      \\
\end{smallmatrix} \right)$, $a\in k_\fp^\times$, for any vector $v\in V_\piv$ in the space of smooth vectors in $\piv$. The image $\CK^\psiv(\piv)$ of this map is called the (smooth) Kirillov model of $\pi_\fp$.

We now describe the structure of $\CK^\psiv(\piv)$ where $k_\fp$ is a non-archimedian local field. 
Let $ \CS(k_\fp^\times)$ be the space of Schwartz functions on $k_\fp^\times$ (i.e.,
locally constant functions of compact support on $k_\fp^\times$). For a supercuspidal representation $\piv$, we have $\CK^\psiv(\piv)=\CS(k_\fp^\times)$. For induced representations $\pi_\fp(\chi_1,\chi_2)$, the space $\CK^\psiv(\piv)$ is a linear span   of $\CS(k_\fp^\times)$ and of a finite-dimensional space $V(\chi_1,\chi_2)$ with $\dim V(\chi_1,\chi_2)= 1$ or $2$. One can take as a basis of $V(\chi_1,\chi_2)$  functions on $k^\times_\fp$ with the support in $\CO_\fp\cap k_\fp^\times$.  More precisely, for an irreducible $\pi_\fp(\chi_1,\chi_2)$, $V(\chi_1,\chi_2)= \bc$-${\rm span}(f_1,f_2)$ with $f_i=\chi_i(a)|a|_\fp^{\haf}\chi_{\CO_\fp}$ and $\chi_{\CO_\fp}$ which is the characteristic function of $\CO_\fp$,  for $\pi_\fp(\chi_1,\chi_2)$ with $\chi_1=\chi_2$, $V(\chi_1,\chi_2)= \bc$-${\rm span}(f_1)$ with 
$f_1=\chi_1(a)|a|_\fp^{\haf}\chi_{\CO_\fp}$, and for for $\pi_\fp(\chi_1,\chi_2)$ with $\chi_1=\chi_2|\cdot|_\fp$, $V(\chi_1,\chi_2)= \bc$-${\rm span}(f_1, f_2)$ with 
$f_1=\chi_1(a)|a|_\fp^{\haf}\chi_{\CO_\fp}$ and $f_2=\chi_2(a)|a|_\fp^{\haf}\log|a|_\fp\chi_{\CO_\fp}$

The action of $\GL(2,k_\fp)$ on $\CK^\psiv(\piv)$ can be described as follows. The action of the Borel subgroup does not depend on the representation (however,  the space of smooth vectors does!), but only on its central character, and is given by
\begin{align}&\piv(n(x))f(a)=\psiv(ax)f(a),\nonumber\\
	&\piv(\bar t)f(a)=f(ta),\\ &\piv(\bar z)f(a)=\om_\piv(z)f(a),\nonumber
\end{align} where $a\in k^\times$,  $\bar t=
\left(\begin{smallmatrix}
t & \\
&    1      \\
\end{smallmatrix} \right)$,  $\bar z=
\left(\begin{smallmatrix}
z & \\
&    z      \\
\end{smallmatrix} \right)$ and $f\in\CK^\psiv(\piv)$. Hence the Whittaker functional $\CW^\psiv$ we started with is given by the evaluation at $a=1$ (i.e., is given by the delta function $\dl_1(f)=f(1)$).

The action of $w$ defines the action of $G$ (via the Bruhat decomposition), and it is known that $\piv(w)$ is given by the integral transform 
\begin{align}\label{Kir-w-act}\piv(w)f(a)=\int_{k_\fp^\times}\om_\piv\inv(t)j_\piv(ta)f(t) d^\times t\ ,\end{align} where the kernel $j_\piv=j_{\piv,\psiv}$ is called the Bessel function of the representation $\piv$. The function $j_\piv$ is a smooth function (i.e., a locally constant for non-archimedian $k_\fp$ and smooth for $k_\fp$ archimedian). We will need a non-trivial bound on $j_\piv$ near $0$. We have
\begin{align}\label{j-bound} |j_\piv(a)|\leq C_\piv|\om_\piv(a)|^{-1/2}|a|_\fp^{-1/4}
\end{align} for $|a|_\fp\leq 1$. This is proved in \cite{Ba}, Corollary 4.2 (see also \cite{JY} for the crucial computation of the germ of the corresponding orbital integral).

\subsection{Kirillov model for $\GLR$}\label{real-Kir-apx}
We recall here the structure of the Kirillov model for unitary representations of $\GLR$. The results we quote are discussed  at length (and proved) in \cite{BM-real} from where we borrow notations as well. We will cover only representations with the {\it trivial} central character. 

Let $\eta\in\{0, 1\}$ and $s\in\bc$, $Re(s)\geq 0$. Let $\Pi_{\eta,s}$ be the (induced)  representation of $G$ in the space of smooth functions $f:G\to\bc$ satisfying $f(n(x)\bar a z(b)h)=sign^\eta(a)|a|^{1/2+s}f(h)$. For $s\not= d-\haf$ where $d$ is a positive integer, the representation $\Pi_{\eta,s}$ is irreducible and we denote it by $\pi_{\eta,s}$.  In that case, it is a unitarazable representation for $Re(s)=0$ (the principal series representations) and for real $s$ satisfying $0<s<1/2$ (the complementary series representations). For $s=d-\haf$ with $d\in\bn$, the representation $\Pi_{\eta,s}$ has the unique irreducible subspace which we denote by $\pi_d$ (suppressing $\eta$ since $\pi_{1,d-1/2}\simeq \pi_{0,d-1/2}$).The quotient space $W_d=V_{\Pi_{1,d-1/2}}/V_{\pi_{1,d-1/2}}$ is finite-dimensional of dimension $2d-1$ (e.g., $W_1\simeq\bc$). The lowest weight vector in the representation $\pi_d$ has the weight $2d$. (Note that in notations of \cite{BM-real}, Theorem 2.5.3, we have $k=2d$.)

Let $\psi_\8(x)=e^{2\pi ix}$. The Bessel function $j_{\pi_d,\psi_\8}$ for representation of discrete series $\pi_d$, $d\in\bn$, is  given by 
\begin{align}\label{j-disc} j_{\pi_d,\psi_\8}(a)=(-1)^d2\pi|a|^\haf J_{2d-1}(4\pi|a|^\haf)\end{align}  for $a>0$ and $j_{\pi_d,\psi_\8}(a)=0$ for $a<0$. Here $J_n$ is the classical $J$-Bessel function (see \cite{Ma}).

For a principal series representation $\pi_{0,ir}$, $ir\in i\br$, we have 
\[j_{\pi_{0,ir},\psi_\8}(x) = -\pi|x|^\haf
\sin(\pi ir)\inv (J_{2ir}(4\pi|x|^\haf) - J_{-2ir}(4\pi|x|^\haf))\] for $x > 0$, and $j_{\pi_{0,ir},\psi_\8}(x) = -\pi|x|^\haf
\sin(\pi ir)\inv (I_{2ir}(4\pi|x|^\haf) - I_{-2ir}(4\pi|x|^\haf))$ for $x<0$. 
Analogous formulas are known for representations of complimentary series (e.g., see \cite{BM-real}). 

We note that for the classical Bessel function one has a well-developed  asymptotic expansion (see \cite{Ma}). In particular, the classical $J$-Bessel function satisfies 
for $ 0 < z \ll  1$, 
\begin{align}\label{j-asymp}
J_\alpha(z)=\sum_{m=0}^N \frac{(-1)^m}{m! \, \Gamma(m+\alpha+1)} {\left(\frac{z}{2}\right)}^{2m+\alpha}+O(|z|^{2N+\al+1})\ ,
\end{align}
 and $J_\alpha(z)\sim\sqrt{\frac{2}{\pi z}}\cos \left(z-\frac{\alpha\pi}{2}-\frac{\pi}{4}\right)$ for $|z|\to\8$ and $Im(z)$ bounded. 
As a result, we have a similar asymptotic expansion at $0$ and $\8$ for  Bessel functions of  all representations.



\begin{thebibliography}{BBB}

\bibitem[Ban]{Ban} W. Banks, A Corollary to Bernstein's Theorem
and Whittaker Functionals on the Metaplectic Group, Math. R. Let. 5 (1998), 781--790. 

\bibitem[Ba]{Ba} M. Baruch,  On Bessel distributions for $\GL(2)$ over a $p$-adic field, J. Number Theory 67 (1997), no. 2, 190--202. 

\bibitem[BM1]{BM} M. Baruch, Z. Mao, Bessel identities in the Waldspurger correspondence over a p-adic field, Amer. J. Math. 125 (2003), no. 2, 225--288.

\bibitem[BM2]{BM-real} M. Baruch, Z. Mao, Bessel identities in the Waldspurger correspondence over the real numbers, Israel J.  Math. 145 (2005), 1-82. 

\bibitem[BS]{BS}  M. Baruch, K. Snitz,  A note on Bessel functions for supercuspidal representations of $\GL(2)$ over a $p$-adic field, Algebra Colloq. 18 (2011), Special Issue No.1, 733--738. 

\bibitem[BZ]{BZ} J. Bernstein, A.V. Zelevinsky, Representations of the group $GL(n,F)$, where $F$ is a local non-Archimedean field, Uspekhi Mat. Nauk 10  (1976), No.3, 5--70. 

\bibitem[Be]{Be} J. Bernstein, Letter to I. Piatetski-Shapiro (1985). In S. Gelbart, I. Piatetski-Shapiro, S. Rallis, Explicit Constructions of Automorphic $L$-Functions, SLN 1254 (1987). 

\bibitem[Bu]{Bu} D. Bump,  Automorphic Forms and Representations, Cambridge University Press, 1998.

\bibitem[E]{E} T. Estermann, On certain functions represented by Dirichlet series, Proc. London Math.
Soc. 27 (1928), 435--448.




\bibitem[G1]{G1} I.~Gelfand, G.~Shilov, Generalized Functions, vol. 1, Academic Press, 1964.



\bibitem[Gr]{Gr} B. Gross, Some applications of Gel'fand pairs to
number theory, Bull. Amer. Math. Soc. (N.S.) 24 (1991), no. 2,
277--301.
\bibitem[GP]{GP} B. Gross and D. Prasad, On the decomposition of a representation of $SO_n$ when restricted to $SO_{n-1}$,
Canad. J. Math. 44 (1992), no. 5, 974--1002.

\bibitem[He]{He} E. Hecke,  \"Uber Modulfunktionen und die Dirichletschen Reihen mit Eulerscher Produktentwicklung, 
Mathematische Annalen 114 (1937), no. 1,  1--28. 

\bibitem[II]{II} A. Ichino, T. Ikeda, On the periods of automorphic forms on special orthogonal groups and the Gross--Prasad conjecture, Geom. Funct. Anal. 19 (2010), no. 5, 1378--1425.

\bibitem[J1]{Jacq-RTF} H. Jacquet, A guide to the relative trace formula, Automorphic Representations, $L$-functions and Applications: Progress and Prospects, 257--272,  Ohio University press, 2005.

\bibitem[J2]{Jacq01} H. Jacquet, Factorization of period integrals, J. Number Theory 87 (2001), no. 1, 109--143.
\bibitem[JL]{JL} H. Jacquet, R. Langlands,  Automorphic forms on $\GLt$, Lecture Notes in Mathematics 114, Berlin, New York: Springer-Verlag
(1970).

\bibitem[JY]{JY} H. Jacquet, Y. Ye, Distinguished representations and quadratic base change for $\GL(3)$, Proc. Amer. Math. Soc. 348(1996), no. 3, 913--939.


\bibitem[K1]{K1} N. 
Kurokawa, 
On the meromorphy of Euler products. I, 
Proc. London Math. Soc. (3) 53 (1986), no. 1, 1--47. 

\bibitem[K2]{K2} N. 
Kurokawa, 
On the meromorphy of Euler products. II, 
Proc. London Math. Soc. (3) 53 (1986), no. 2, 209--236. 

\bibitem[LM1]{LM1} E. Lapid, Z. Mao, On Whittaker-Fourier coefficients of automorphic forms on $\tilde{Sp}_n$,  J. Number Theory 146 (2015), 448–505.

\bibitem[LM2]{LM2} E. Lapid, Z. Mao, Model transition for representations of metaplectic type. With an appendix by M. Tadi\'{c}.  IMRN 2015, no. 19, 9486--9568.

\bibitem[LM3]{LM3} E. Lapid, Z. Mao, Stability of certain oscillatory integrals. Int. Math. Res. Not. IMRN 2013, no. 3, 525--547.

\bibitem[Ma]{Ma} W. Magnus et al., Formulas and Theorems for the Special Functions, Springer, 1966.

\bibitem[MZ]{MZ} D. Montgomery, L. Zippin, Topological Transformation Groups, Interscience Publishers, 1955. 

\bibitem[PS]{PS} I. Piatetsky-Shapiro, Euler subgroups.
Lie groups and their representations, Halsted, (1975), 597--620.

\bibitem[Po]{Po} A. Popa, Whittaker newforms for archimedian representations of $\GLt$, J.  Number Theory, 128 (2008), 1637--1645. 

\bibitem[Ra]{Ra} S. Ramanujan,  On certain arithmetical functions. Trans. Cambridge Philos. Soc. 22 (9)(1916), 159--184.

\bibitem[SV]{SV}  Y. Sakellaridis, A. Venkatesh, 
Periods and harmonic analysis on spherical varieties,  Asterisque, 396 (2017).


\bibitem[S]{S} D. Soudry, The $L$ and $\g$ factors for generic representations of $GSp(4,k)\times GL(2,k)$ over a local non-archimedian field $k$, Duke Math. J. 51 (1984), 355--394. 

\bibitem[Ta]{Ta} T. Tamagawa, Ad\`{e}les, Algebraic Groups and Discontinuous Subgroups, Proc. Sympos. Pure Math. IX, Providence, R.I.: American Mathematical Society, 113--121 (1966).

\bibitem[T]{T} T. Tao, Hilbert's fifth problem and related topics,  American Mathematical Society. 2014.

\bibitem[Wa]{Wa}  J.-L.Waldspurger, Sur les valeurs de fonctions $L$-automorphes en leur centre de sym\'{e}trie, Comp. Math. t. 54 (1985),173--242.



\bibitem[We]{We} A. Weil, Adeles and Algebraic Groups, Progress in Mathematics 23, Boston, MA: Birkhäuser Boston, (1961).

\bibitem[Za]{Za} D. Zagier, Modular forms whose Fourier coefficients involve zeta-functions of quadratic fields. Modular functions of one variable, VI , 105--169. Lecture Notes in Math., Vol. 627, Springer, Berlin, (1977).


\end{thebibliography}
\end{document}